\documentclass[12pt]{article}
\usepackage{e-jc+}

\title{Quantitative Combinatorial Nullstellensatz}

\author{Uwe Schauz\\
\small Department of Mathematical Science\\[-0.8ex]
\small Xi’an Jiaotong-Liverpool University\\[-0.8ex]
\small Suzhou 215123, China\\[-0.8ex]
\small \texttt{uwe.schauz@xjtlu.edu.cn}}

\date{Nov 14, 2006\\
\small Mathematics Subject Classifications:
  41A05, 13P10, 05E99, 11C08, 11D79, 05C15, 15A15}

\usepackage{fixltx2e}

\usepackage[ansinew]{inputenc}

\usepackage{datetime}
\usepackage{ngerman}
\usepackage{eurosym}
\usepackage{textcomp}

\usepackage{epic}
\usepackage{eepic}

\usepackage{amsmath}
\usepackage{amsthm}

\usepackage{amssymb}
\usepackage{stmaryrd}

\usepackage{enumerate}

\usepackage{afterpage}
\usepackage{xspace}

\newtheoremstyle{Theorem}{.7\baselineskip}{1\baselineskip}{\itshape}{}{\bfseries}{.}{ }{}
\theoremstyle{Theorem}
\newtheorem{Satz}{Theorem}[section]
\newtheorem{EquiDef}[Satz]{Equivalence and Definition}
\newtheorem{Korollar}[Satz]{Corollary}

\newtheorem{Proposition}[Satz]{Proposition}

\newtheorem{Lemma}[Satz]{Lemma}
\newtheorem{Vermut}[Satz]{Conjecture}

\newtheoremstyle{Definition}{.6\baselineskip}{.8\baselineskip}{}{}{\bfseries}{.}{ }{}
\theoremstyle{Definition}
\newtheorem{Definition}[Satz]{Definition}
\theoremstyle{definition}
\newtheorem{Beispiel}[Satz]{Example}

\theoremstyle{remark}
\newtheorem{Bemerkung}[Satz]{Remark}

\newenvironment{Beweis}[1][Proof]{\begin{proof}[#1]}{\end{proof}}
\newcommand\ps{\small}
\newenvironment{proofsize}[1]{\begin{small}#1}{\end{small}}

\newenvironment{Liste}[1][]
  {\begin{list}{#1}{\setlength{\topsep}{\itemsep}
   \setlength{\leftmargin}{.5em}\setlength{\labelwidth}{.5em}\setlength{\labelsep}{0em}
   \setlength{\listparindent}{\parindent}\setlength{\parsep}{0ex}}}
  {\end{list}}

  {\setcounter{enumii}{0}\begin{Liste}[\stepcounter{enumii}%
   \emph{#1\arabic{enumi}.\arabic{enumii}}]}
  {\end{Liste}}

\newcommand\Index[1]{\index{#1}#1}
\newcommand\dashed[1]{\hbox{-- }#1\hbox{ --}\xspace}
\newcommand\´{\nolinebreak\hspace{.08em}}
\newcommand\ms{\hspace{\mathsurround}}
\newcommand\noms{\hspace{-\mathsurround}}
\newcommand\vsp{\vspace{1ex}}
\newcommand\vspp{\vspace{3ex}}

\newcommand\Rand[1]{
  \marginpar{\raggedleft\scriptsize\hspace{0pt}#1}}%

\renewcommand{\(}{$\noms}
\renewcommand{\)}{\noms$}
\renewcommand\frac[2]{\genfrac{}{}{.4pt}{}{#1}{#2}}
\renewcommand\tfrac[2]{\genfrac{}{}{.4pt}{1}{#1}{#2}}

\newcommand\Strut{\phantom{\rlap{\ensuremath{\displaystyle\sum}}}}
\newcommand\mathRand[1]{\hspace{\mathsurround}\Rand{#1}\nolinebreak\noms}
\def\rand #1"#2"{\mathRand{\(#2\)}#1#2}
\def\randd #1"#2"#3\randd#4"#5"{\mathRand{\(#2\), \(#5\)}#1#2#3#4#5}

\newcommand\eqby[2][=]%
  {\ensuremath{\overset{\makebox[0pt]{\ensuremath{\smash[t]{\scriptstyle#2}}}}{#1}}}
\newcommand\Ceqby[2][=]{\quad\eqby[#1]{#2}\quad}
\newcommand\ceqby[2][=]{&\Ceqby[#1]{#2}}
\newcommand\Ceq[1][=]{\quad#1\quad}
\newcommand\ceq[1][=]{&\Ceq[#1]}

\renewcommand\o{\varnothing}

\newcommand\F[1][]{\mathbb{F}_{\!#1}}
\newcommand\C{\mathbb{C}}
\newcommand\R{\mathcal{R}}
\newcommand\Rl{\mathbb{R}}

\newcommand\Z{\mathbb{Z}}
\newcommand\N{\mathbb{N}}

\renewcommand\sb{\subseteq}
\renewcommand\sp{\supseteq}
\newcommand\nin{\notin}

\newcommand\sm{\setminus}
\newcommand\ssm{\textup{\texttt{\char92}}}
\newcommand\mto{\mapsto}
\newcommand\lmto{\longmapsto}
\newcommand\lto{\longrightarrow}
\newcommand\To{\Rightarrow}
\newcommand\lTo{\Longrightarrow}
\newcommand\Eqi{\Leftrightarrow}
\newcommand\lEqi{\Longleftrightarrow}

\newcommand\fa{\forall\,}
\newcommand\ex{\exists\,}

\newcommand\DP{\colon\discretionary{\!\kern -.17em}{}{}}
\newcommand\mitsymbol{\textup{\textbrokenbar}}
\renewcommand\mit{\,\ \discretionary{\mitsymbol}{}{}\mitsymbol\ \,}
\renewcommand\div{\mathrel{\bigm\lfloor\!\!\!\bigm\lfloor}}
\newcommand\vid{\mathrel{\bigm\rfloor\!\!\!\bigm\rfloor}}
\newcommand\ndiv{\mathrel{\;\!\div\hspace{-12pt}\kern0pt\lower2pt%
  \hbox{\ensuremath{^\diagup}}\!}}
\newcommand\ndivps{\mathrel{\;\!\div\hspace{-9pt}\kern0pt\lower2pt%
  \hbox{\ensuremath{^\diagup}}\!}}
\newcommand\nvid{\mathrel{\;\!\vid\hspace{-12pt}\kern0pt\lower2pt%
  \hbox{\ensuremath{^\diagup}}\!}}
\newcommand\nvidps{\mathrel{\;\!\vid\hspace{-9pt}\kern0pt\lower2pt%
  \hbox{\ensuremath{^\diagup}}\!}}
\newcommand\Div{\mathrel{\Bigm\lfloor\hspace{-6pt}\Bigm\lfloor}}

\newcommand\nDiv{\mathrel{\,\Div\hspace{-12.13pt}\diagup}}

\newcommand\n[1][n]{{\phantom{|}}^{\!\!\!\!#1}}

\providecommand\abs[1]{\lvert#1\rvert}
\providecommand\Abs[1]{\bigl\lvert#1\bigr\rvert}

\newcommand\?[1]{{?}_{(#1)}}

\DeclareMathOperator\supp{supp}
\DeclareMathOperator\Bild{Im}

\DeclareMathOperator\Pot{\mathcal{P}}

\DeclareMathOperator\per{per}
\DeclareMathOperator\Adj{Adj}

\newcommand\X{\mathfrak{X}}
\newcommand\PX[1][\X]{P\!/\´\!#1}

\newcommand\x{\chi}
\renewcommand\S{\mathcal{S}}
\newcommand\St{\S_{\textup{triv}}}
\newcommand\D{\mathcal{D}}

\newcommand\p{\varphi}

\newcommand\tx{\skew{2}{\tilde}{x}}
\newcommand\hx{\skew{2}{\hat}{x}}
\renewcommand\d{\delta}

\newcommand\e{\varepsilon}

\newcommand\Chi{\kern0pt\lower-2.5pt\hbox{\ensuremath\chi\!}}

\newcommand\E{\mathsf{1}}

\newcommand\sel[1]{\langle#1\rangle}
\newcommand\sTo{{{\shortrightarrow}\!\!\!\!\!\hspace{.75pt}%
  {\shortrightarrow}\!\!\!\!\!\hspace{.75pt}{\shortrightarrow}\!\!\!\!\!\hspace{.75pt}%
  {\shortrightarrow}\!\!\!\!\!\hspace{.75pt}{\shortrightarrow}}}
\newcommand\sFrom{{{\shortleftarrow}\!\!\!\!\!\hspace{.75pt}%
  {\shortleftarrow}\!\!\!\!\!\hspace{.75pt}{\shortleftarrow}\!\!\!\!\!\hspace{.75pt}%
  {\shortleftarrow}\!\!\!\!\!\hspace{.75pt}{\shortleftarrow}}}

\newcommand\sto{\scriptscriptstyle{\sTo}}
\newcommand\sfrom{\scriptscriptstyle{\sFrom}}

\newcommand\G[1]{{\overset{\makebox(0,-1.9){\ensuremath{#1\hspace{-.2em}}}}{G}}}

\newcommand\ed[1]{e^{\!#1}}

\begin{document}
\maketitle

\begin{abstract}
The main result of this paper is a coefficient formula that sharpens
and generalizes Alon and Tarsi's Combinatorial Nullstellensatz. On
its own, it is a result about polynomials, providing some
information about the polynomial map
$P|_{\X_1\times\dotsb\times\X_n}$ when only incomplete information
about the polynomial $P(X_1,\dotsc,X_n)$ is given.

In a very general working frame, the grid points
$x\in\X_1\times\dotsb\times\X_n$ \linebreak which do not vanish
under an algebraic solution -- a certain describing polynomial
$P(X_1,\dotsc,X_n)$ -- \,correspond to the explicit solutions of a
problem. As a consequence of the coefficient formula, we prove that
the existence of an algebraic solution is equivalent to the
existence of a nontrivial solution to a problem. By a problem, we
mean everything that ``owns'' both, a set $\S$, which may be called
the \emph{set of solutions}; and a subset $\St\subseteq\S$, the
\emph{set of trivial solutions}.

We give several examples of how to find algebraic solutions, and how
to apply our coefficient formula. These examples are mainly from
graph theory and combinatorial number theory, but we also prove
several versions of Chevalley and Warning's Theorem, including a
generalization of Olson's Theorem, as examples and useful
corollaries.

We obtain a permanent formula by applying our coefficient formula to
the matrix polynomial, which is a generalization of the graph
polynomial. This formula is an integrative generalization and
sharpening of\´:

1. Ryser's permanent formula.

2. Alon's Permanent Lemma.

3. Alon and Tarsi's Theorem about orientations and colorings of
graphs.\\[2pt]
Furthermore, in combination with the Vigneron"=Ellingham"=Goddyn
property of planar \(n\)"~regular graphs, the formula contains as
very special cases\´:

4. Scheim's formula for the number of edge \(n\)-colorings of such
graphs.

5. Ellingham and Goddyn's partial answer to the list coloring
conjecture.
\end{abstract}

\newpage
\section*{Introduction}
Interpolation polynomials $P=\sum_{\d\in\N^n}P_\d X^\d$ on finite
``grids''  \rand$"\X":=\X_1\times\dotsb\times\X_n\,\sb\,\F\,\,\n$
are not uniquely determined by the interpolated maps
\rand$"P|_\X"\DP x\mto P(x)$. One could restrict the partial degrees
to force the uniqueness. If we only restrict the total degree to
\linebreak[4]$\deg(P)\leq d_1+\dotsb+d_n$, where
\rand$"d_j":=\abs{\X_j}-1$, the interpolation polynomials $P$ are
still not uniquely determined, but they are partially unique. That
is to say, there is one (\´and in general only one\´) coefficient in
$P=\sum_{\d\in\N^n}P_\d X^\d$ that is uniquely determined, namely
\rand$"P_d"$ with $d:=(d_1,\dotsc,d_n)$. We prove this in
Theorem\,\ref{sz2.cn} by giving a formula for this coefficient. Our
coefficient formula contains Alon and Tarsi's Combinatorial
Nullstellensatz\,\cite[Th.\,1.2]{al2}, \cite{al3}\´:
\begin{equation}\label{eq.cnit}
P_d\,\neq\,0\quad\lTo\quad P|_\X\,\not\equiv\,0\,\ .
\end{equation}\smallskip

This insignificant"=looking result, along with Theorem\,\ref{sz2.cn}
and its corollaries \ref{kor.cn}, \ref{kor2.cn} and \ref{kor4.cn},
are astonishingly flexible in application. In most applications, we
want to prove the existence of a point $x\in\X$ such that
$P(x)\neq0$. Such a point $x$ then may represent a coloring, a graph
or a geometric or number"=theoretic object with special properties.
In the simplest case we will have the following correspondence\´:
\begin{equation}\label{sch.sol}
\begin{split}
  \text{$\X$}\quad
    &\longleftrightarrow\quad\text{Class of Objects}\\
  \text{$x$}\´\quad
    &\longleftrightarrow\quad\text{Object}\\
  P(x)\neq0\,\quad
    &\longleftrightarrow\quad\text{``Object is interesting (\´a solution\´).''}\\
  P|_\X\not\equiv0\,\quad
    &\longleftrightarrow\quad\text{``There exists an interesting object
                                   (\´a solution\´).''}
\end{split}
\end{equation}\smallskip

\noindent This explains why we are interested in the connection
between $P$ and $P|_\X$: In general, we try to retrieve information
about the polynomial map $P|_\X$ using incomplete information about
$P$. One important possibility is if there is (\´exactly\´) one
trivial solution $x_0$ to a problem, so that we have the information
that $P(x_0)\neq0$. If, in this situation, we further know that
$\deg(P)<d_1+\dotsc+d_n$, then Corollary\,\ref{kor.cn} already
assures us that there is a second (\´nontrivial\´\´) solution
$x\!$\´, i.e., an $x\neq x_0$ in $\X$ such that $P(x)\neq0$. The
other important possibility is that we do not have any trivial
solutions at all, but we know that $P_d\neq0$ and $deg(P)\leq
d_1+\dotsc+d_n$. In this case, $P|_\X\not\equiv0$ follows from
\eqref{eq.cnit} above or from our main result,
Theorem\,\ref{sz2.cn}\,. In other cases, we may instead apply
Theorem\,\ref{sz.cn}\,, which is based on the more general concept
from Definition\,\ref{def.dl} of \(d\)"~leading coefficients.

In Section\,\ref{sec.fap}\´, we demonstrate how most examples from
\cite{al2} follow easily from our coefficient formula and its
corollaries.
The new, quantitative version\,\ref{sz2.cn}\,\((i)\) of the
Combinatorial Nullstellensatz is, for example, used in
Section\,\ref{sec.mp}, where we apply it to the matrix polynomial --
a generalization of the graph polynomial -- to obtain a permanent
formula. This formula is a generalization and sharpening of several
known results about permanents and graph colorings (\´see the five
points in the abstract\´). We briefly describe how these results are
derived from our permanent formula.

We show in Theorem\,\ref{sz.as} that it is theoretically always
possible, both, to represent the solutions of a given problem $\Pot$
(\´see Definition\,\ref{def.prob}\´) through some elements $x$ in
some grid $\X\!$, and to find a polynomial $P\!$, with certain
properties (\´e.g., $P_d\neq0$ as in \eqref{eq.cnit} above\´), that
describes the problem\´:
\begin{equation}
P(x)\neq0\quad\lEqi\quad\text{``$x$ represents a solution of
$\Pot\!$.''}
\end{equation}
We call such a polynomial $P$ an \emph{algebraic solution} of
$\Pot\!$, as its existence guarantees the existence of a nontrivial
solution to the problem $\Pot$.

Sections \ref{sec.fap} and \ref{sec.mp} contain several examples of
algebraic solutions. Algebraic solutions are particularly easy to
find if the problems possess exactly one trivial solution: due to
Corollary\,\ref{kor.cn}\´, we just have to find a describing
polynomial $P$ with degree \linebreak[4]$\deg(P)<d_1+\dotsc+d_n$ in
this case. Loosely speaking, Corollary\,\ref{kor.cn} guarantees that
every problem which is not too complex, in the sense that it does
not require too many multiplications in the construction of $P\!$,
does not possess exactly one (\´the trivial\´) solution.

In Section\,\ref{sec.cn} we give a slight generalization of the
(\´first\´) Combinatorial Nullstellensatz -- a sharpened
specialization of Hilbert's Nullstellensatz -- and a discussion of
Alon's original proving techniques. Note that, in
Section\,\ref{sec.cf} we used an approach different from Alon's to
verify our main result. However, we will show that Alon and Tarsi's
so"=called polynomial method can easily be combined with
interpolation formulas, such as our inversion formula\,\ref{sz.itg},
to reach this goal.

Section\,\ref{sec.ZZm} contains further generalizations and results
over the integers $\Z$ and over $\Z/m\Z$. Corollary\,\ref{kor3.cn}
is a surprising relative to the important Corollary\,\ref{kor.cn}\´,
one which works without any degree restrictions.
Theorem\,\ref{kor4.cn}\´, a version of Corollary\,\ref{kor2.cn}\´,
is a generalization of Olson's Theorem.

Most of our results hold over integral domains, though this condition has been weakened in
this paper (\´see \ref{eqdef.itg} for the definition of integral grids\´). In the important case of
the Boolean grid $\X=\{0,1\}^n$\!\!, our results hold over arbitrary commutative rings $\R$.
Our coefficient formulas are based on the interpolation formulas in Section\,\ref{sec.ip}\,,
where we generalize known expressions for interpolation polynomials over fields to
commutative rings $\R$. We frequently use the constants and definitions from
Section\,\ref{sec.nc}\,.

For newcomers to this field, it might be a good idea to start with
Section\,\ref{sec.fap} to get a first impression.

We will publish two further articles\´: One about a sharpening of
Warning's classical result about the number of simultaneous zeros of
systems of polynomial equations over finite fields \cite{sch2}, the
other about the numerical aspects of using algebraic solutions to
find explicit solutions, where we present two polynomial"=time
algorithms that find nonzeros of polynomials \cite{sch3}.

This paper was first published in the Electronic Journal of Combinatorics\,\cite{schAlg} under
another title, please use this reference for citations.

\clearpage

\section{Notation and constants}\label{sec.nc}

\index{$\R$, $\Z_m$, $\F[p^k]$, $\N$}%
\index{$n\vid p$, $(n]$, $[n)$, $[n]$}%
\index{$\?{\mathcal{A}}$, $\Pi$, $\Sigma$, $\supp(y)$, $\bigotimes$}%
  \rand$\noms"\R"$ is always a 
  commutative ring with $1\neq0$.\\ \randd$"\F[p^k]"$ denotes the field with $p^k$
  elements (\´$p$ prime\´) and \randd$"\Z_m":=\Z/m\Z$.\\
  We write \randd$"p\div n"$ (\´or $n\vid p$\´)
  for ``$p$ \emph{divides} $n$'' and abbreviate
  \randd$"\S\ssm s":=\S\sm\{s\}$\medskip.\\
  For \hspace{-.7pt}\rand$n\in"\N":=\{0,1,2,\dotsc\}$
  we set\´:\\
 {\newcommand\pl{\makebox[3pt]{(}}%
  \newcommand\pr{\makebox[3pt]{)}}%
  \newcommand\bl{\makebox[3pt]{[}}%
  \newcommand\br{\makebox[3pt]{]}}%
  \rand$"\pl n\br"=\pl0,n\br:=\{1,2,\dotsc,n\}$,\\
  \rand$"\bl n\pr"=\bl0,n\pr:=\{0,1,\dotsc,n\!-\!1\}$,\\
  \rand$"\bl n\br"=\bl0,n\br:=\{0,1,\dotsc,n\}$. (\´Note that $0\in[n]$.\´)\bigskip
 }

  \noindent
  For statements $\mathcal{A}$ the ``Kronecker query''
  $\?{\mathcal{A}}$ is defined by\´:\\
  \rand$"\?{\mathcal{A}}":=\begin{cases}
    0&\text{if $\mathcal{A}$ is false,}\\
    1&\text{if $\mathcal{A}$ is true.}
  \end{cases}$\\
\medskip

  \noindent
  For finite tuples (\´and maps\´) $d=(d_j)_{j\in J}$ 
  and sets $\Gamma $ we define\´:\\
  \randd$"\Pi d"
  :=\prod_{j\in J}d_j$ , \ \randd$"\Pi\Gamma" :=\prod_{\gamma\in\Gamma}\gamma$
  \ and\\
  \randd$"\Sigma d"
  :=\sum_{j\in J}d_j$ ,
  \ \randd$"\Sigma\Gamma" :=\sum_{\gamma\in\Gamma}\gamma$ .\bigskip

  For maps $y,z\DP\X\lto\R$ with finite domain we identify the map
  \rand$"y"\DP x\lmto y(x)$ with \smallskip the tuple $(y(x))_{x\in\X}\,\in\,\R^\X$\!.
  Consequently, the product with matrices
  $\Psi=(\psi_{\d,x})\,\in\,\R^{D\times\X}$
  is given by
  \rand$"\Psi y":=\bigl(\´\sum_{x\in\X}\psi_{\d,x}\,y(x)\,\bigr)_{\d\in D}\in\R^D$%
  \!\smallskip.\\
  We write \randd$"yz"$ for the pointwise product, $(yz)(x):=y(x)z(x)$.
  If nothing else is said, \randd$"y^{-1}"$ is also defined pointwise,
  $y^{-1}(x):=y(x)^{-1}\!$, if $y(x)$ is invertible for all $x\in\X$.\smallskip\\
  We define \rand$"\supp(y)":=\{\,x\in\X\mit y(x)\neq0\,\}$.
  \bigskip

  The tensor product \rand$"\bigotimes"_{j\in(n]}y_j$ of maps $y_j\DP\X_j\lto\R$
  is a map \smallskip $\X_1\times\dotsb\times\X_n\lto\R$, it is defined by
  $(\bigotimes_{j\in(n]}y_j)(x):=\prod_{j\in(n]}y_j(x_j)$\medskip.\\
  Hence, the tensor product $\bigotimes_{j\in(n]}a^j$ of tuples
  $a^j:=(a^j_{x_j})_{x_j\in\X_j}$, $j\in(n]$,
  is the tuple
  $\bigotimes_{j\in(n]}a^j\,:=\,
  \bigl(\´\prod_{j\in(n]}a^j_{x_j}\bigr)_{x\in\X_1\times\dotsb\times\X_n}$\medskip.\\
  The tensor product $\bigotimes_{j\in(n]}\Psi^j$ of matrices
  $\Psi^j=(\psi^j_{\d_j,x_j})\smash{_{\d_j\in D_j\atop x_j\in\X_j}}$, $j\in(n]$,
  is the matrix
  $\bigotimes_{j\in(n]}\Psi^j\,:=\,
  \bigl(\prod_{j\in(n]}\psi^j_{\d_j,x_j}\bigr)%
  _{\d\in D_1\times\dotsb\times D_n\atop x\in\X_1\times\dotsb\times\X_n}$\medskip.\\
  Tensor product and matrix"=tuple multiplication go well together\´:
\begin{multline}\label{tens}
  \!\!\!\!\bigl(\bigotimes_{j\in(n]}\Psi^j\,\bigr)\bigotimes_{j\in(n]}a^j
  \,=\,\bigl(\prod_{j\in(n]}\psi^j_{\d_j,x_j}\bigr)_{\!{\d\in D\atop x\in\X\,\,}}
   \bigl(\prod_{j\in(n]}a^j_{x_j}\bigr)_{x\in\X}
  \,=\,\Bigl(\´\sum_{x\in\X}\prod_{j\in(n]}\psi^j_{\d_j,x_j}a^j_{x_j}\Bigr)_{\d\in D}\\
  \,=\,\Bigl(\prod_{j\in(n]}\sum_{x_j\in\X_j}\psi^j_{\d_j,x_j}a^j_{x_j}\Bigr)_{\d\in D}
  \,=\,\bigotimes_{j\in(n]}%
       \bigl(\sum_{x_j\in\X_j}\psi^j_{\d_j,x_j}a^j_{x_j}\bigr)_{\d_j\in D_j}
  \,=\,\bigotimes_{j\in(n]}(\Psi^j a^j)\ .
\end{multline}\smallskip\pagebreak[2]

In the whole paper we work over Cartesian products
$\X:=\X_1\times\dotsb\times\X_n$ of subsets $\X_j\sb\R$ of size
$d_j+1:=\abs{\X_j}<\infty$. We define\´:\medskip\smallskip

\index{grid}\index{\(d\)-grid}%
\index{$\X$, $[d]$, $d$}
\begin{Definition}[\(d\)"~grids $\X$]\ \vspace{1.2ex}%
\label{def.X}
 \Rand{}
 \Rand{\(\X,[d]\)}
 \Rand{\(d=d(\X)\)}\ \vspace{1.2ex}\\
\noindent\renewcommand\arraystretch{1.8}
\begin{tabular*}{\linewidth}[t]%
  {@{\,}p{.47\linewidth}@{\,\ \vrule width .2pt}@{\ \ }p{.47\linewidth}@{\!}}%
  \raisebox{2ex}{For all $j\in(n]$ we define\´:}
     & \raisebox{2ex}{In n dimensions we define\´:}\\[-2ex]\hline
  $\X_j\sb\R$ \,is always a finite set $\neq\o$.
     & $\X:=\X_1\times\dotsb\times\X_n\,\sb\,\R^n$
       \,is a \emph{\(d\)"~grid} for\\[-.2em]
  \fbox{$d_j=d_j(\X_j):=\abs{\X_j}-1$} \ and
     & $d=d(\X):=(d_1,\dotsc,d_n)$\,.\\[-.2em]
  $[d_j]:=\{0,1,\dotsc,d_j\}$\,.
     & $[d]:=[d_1]\!\times\!\dotsb\!\times\hspace{-1.5pt}[d_n]$
       \ is \!a \(d\)"~grid in $\Z^n$\!.\\
\end{tabular*}
\end{Definition}\vspp

The following function $N\DP\X\lto\R$ will be used throughout the
whole paper. The $\psi_{\d,x}$ are the coefficients of the Lagrange
polynomials $L_{\X,x}$, as we will see in Lemma\,\ref{lem.lagr}\,.
We define\´:\medskip\smallskip

\index{$N$, $\Psi$, $L_{\X,x}$,$e_x$}
\begin{Definition}[$N_\X$, $\Psi_\X$, $L_{\X,x}$ and $e_x$]%
\label{def.N}\ \\[.1ex]
Let $\X:=\X_1\times\dotsb\times\X_n\,\sb\,\R^n$ be a \(d\)"~grid,
i.e.,
$d_j=\abs{\X_j}-1$ for all $j\in(n]$.\\
 \Rand{}
 \Rand{\(e_x\´,\,L_{\X,x}\)}
 \Rand{\(N,\Psi\)}\ \vspace{1.2ex}\\
\noindent\renewcommand\arraystretch{1.8}
\begin{tabular*}{\linewidth}[t]
  {@{\,}p{.47\linewidth}@{\,\ \vrule width .2pt}@{\ \ }p{.47\linewidth}@{\!}}
  \raisebox{2ex}{For $x\in\X_j$ and $\d\in[d_j]$ we set\´:}
     & \raisebox{2ex}{For $x\in\X$ and $\d\in[d]$ we set\´:}\\[-2ex]\hline
  $\,e^j_x\DP\X_j\to\R$\,,\quad$e^j_x(\tx)\´:=\´\?{\tx=x}$ \,.
     & $\,e_x:=\´\bigotimes_{j\in(n]}e^j_{x_j}\,=\´(\´\tx\mto\?{\tx=x}\´)$ \,.\\
  $\,L_{\X_j\!\ssm x}(X)\,:=\,\prod_{\hx\in\X_j\!\ssm x}(X-\hx)$ \,.
     & $\,L_{\X,x}(X_1,\dotsc,X_n)\,:=\,\prod_{j}L_{\X_j\!\ssm x_j}(X_j)$ .\\
  $\,N_j=N_{\X_j}\DP\X_j\lto\R$ \,is defined by\´:
     & $\,N=N_\X\DP\X\lto\R$ \,is defined by\´:\\[-.2em]
  \fbox{\parbox[c]{.85\linewidth}{\strut\´$N_j(x)\,:=\,L_{\X_j\!\ssm x}(x)\ .
        \phantom{\rlap{\ensuremath{\bigotimes_{j\in(n]}}}}$}}
     & \fbox{\parbox[c]{.85\linewidth}{\strut\´$N\,:=\,\bigotimes_{j\in(n]}N_j
        \,=\bigl(\´x\mto L_{\X,x}(x)\´\bigr)\,.\!\!\!\!\!\!$}}\\[1.4em]
  \fbox{\parbox[c]{.85\linewidth}
  {$\newline\smallskip\displaystyle\,\Psi^j%
    \,:=\,(\psi^j_{\d,x})_{\d\in[d_j]\!\atop x\in\X_j\,\,}$\quad
  with
  \newline$\displaystyle\mbox{}\,\psi^j_{\d,x}
  \,:=\,\sum_{{\,\,\,\Gamma\sb\X_j\!\ssm x\atop\abs{\Gamma}=d_j-\d}}\!\!\Pi(-\Gamma)$
  \vspace{2.2ex}
  \newline\mbox{}\,and in particular \,$\displaystyle\psi^j_{d_j,x}=1$\,.}}
  \refstepcounter{equation}\hfill(\arabic{equation})\label{eqP1}
    & \fbox{\parbox[c]{.85\linewidth}
      {$\newline\smallskip\displaystyle\,\Psi
       \,=\,(\psi_{\d,x})_{\d\in[d]\atop x\in\X\,\,}
       \,:=\,\bigotimes_{j\in(n]}\Psi^j
       $\ \,, \,i.e.,
      \newline$\displaystyle\mbox{}\,\psi_{\d,x}
      \,:=\,\prod_{j\in(n]}\psi^j_{\d_j,x_j}$%
      \vspace{2.565ex}
      \newline\mbox{}\,and in particular \,$\displaystyle\psi_{d,x}=1$\,.%
                                         \phantom{$\psi^j_{d_j,x}$}
      }}
      \refstepcounter{equation}\hfill(\arabic{equation})\label{eqP2}\\
\end{tabular*}
\end{Definition}\vspp\pagebreak[2]\smallskip

We use multiindex notation for polynomials, i.e.,
 \rand$"X^{(\d_1,\dotsc,\d_n)}":=X_1^{\d_1}\dotsm X_n^{\d_n}$
and we define \rand$"P_\d=(P)_\d"$ to be the coefficient of $X^\d$
in the standard expansion of
\rand$P\,\in\,"\R[X]":=\R[X_1,\dotsc,X_n]$. That means
$P=P(X)=\sum_{\d\in\N^n}P_\d X^\d$ and
$(X^\e)_\d=\?{\d=\e}$.\pagebreak[2]

Conversely, for tuples $P=(P_\d)_{\d\in\D}\,\in\,\R^\D$,
we set
\rand$"P(X)":=\sum_{\d\in\D}P_\d X^\d$\!.
 In this way we identify the set of tuples
\rand$"\R^{[d]}"=\R^{[d_1]\times\dotsb\times[d_n]}$ with
\rand$"\R[X^{\leq d}]"$, the set of polynomials
 $P=\sum_{\d\leq d}P_\d X^\d$ with restricted partial degrees
$\deg_j(P)\leq d_j$. It will be clear from the context whether we
view $P$ as a tuple $(P_\d)$ in $\R^{[d]}$, a map $[d]\lto\R$ or a
polynomial $P(X)$ in $\R[X^{\leq d}]$. \rand$"P(X)|_\X"$ stands for
the map $\X\lto\R$, $x\lmto P(x)$.\bigskip

\noindent \underline{We have introduced the following four
related or identified objects\´:}%
\index{$P_\d$, $P(X)$, $P"|_\X=P(X)"|_\X$}
\begin{align}
 &\text{Maps\´:}       &&
  \text{Tuples\´:}     &&
  \text{Polynomials\´:}&&
  \text{Polynomial Maps\´:}\notag\\
 &\d\mto P_\d\,,      &\ \ \,\qquad&
  P=(P_\d)            &\ \ \,\qquad&
  P(X)={\textstyle\sum}\,P_\d X^\d &\ \ \,\qquad&
  P(X)|_\X\DP\´x\mto P(x)\,,\qquad\ \notag\\[-2pt]
 &[d]\to\R\!\!\!      &&
  \!\in\R^{[d]}       &&
  \!\in\R[X^{\leq d}] &&
  \X\to\R
\end{align}\

With these definitions we get the following important formula\´:

\index{Lagrange polynomial}
\begin{Lemma}[Lagrange polynomials]\label{lem.lagr}%
$$\boxed{\,(\Psi e_x)(X)
  \,:=\,\sum_{\d\in[d]}\psi_{\d,x}X^\d
  \,=\,\prod_{j\in(n]}^{\null}
     \,\prod_{\hx_j\in\X_j\!\ssm x_j\!\!\!\!\!\!\!\!\!}(X_j-\hx_j)
  \,=:\,L_{\X,x}}\,\ .$$
\end{Lemma}\vspp

\begin{Beweis}We start with the
one"=dimensional case. Assume $x\in\X_j$, then
\begin{equation}\ps\label{1L}
\begin{split}
  (\Psi^j e^j_x)(X_j)
  \ceq\bigl(\sum_{\d\in[d_j]}\psi^j_{\d,x}X_j^\d\,\,\bigr)\\
  \ceq\sum_{\d\in[d_j]}\,
    \sum_{{\,\,\,\Gamma\sb\X_j\!\ssm x\atop\abs{\Gamma}=d_j-\d}}
    \!\!X_j^\d\,\Pi(-\Gamma)\\
  \ceq\sum_{\hat\Gamma\sb\X_j\!\ssm x\!\!\!\!}
    X_j^{\abs{(\X_j\!\ssm x)\sm\hat\Gamma}}\,\Pi(-\hat\Gamma)\\
  \ceq\prod_{\hx\in\X_j\!\ssm x\!\!\!\!}(X_j-\hx)\ .
\end{split}
\end{equation}
In $n$ dimensions and for $x\in\X$ we conclude\´:
\begin{equation}\ps
\begin{split}
  (\Psi e_x)(X)
  \ceq\Bigl(\bigl(\bigotimes\nolimits_{\!j}\Psi^j\,\bigr)
                       \bigotimes\nolimits_{\!j}e^j_{x_j}\,\Bigr)(X)\\
  \ceqby{\eqref{tens}}\Bigl(\bigotimes\nolimits_{\!j}%
                             \bigl(\Psi^j e^j_{x_j}\bigr)\,\Bigr)(X)\\
  \ceq\prod\nolimits_j\bigl((\Psi^j e^j_{x_j})(X_j)\bigr)\\
  \ceqby{\eqref{1L}}\prod_{j\in(n]}\,%
        \prod_{\hx_j\in\X_j\!\ssm x_j\!\!\!\!\!\!\!\!\!}(X_j-\hx_j)\ .
\end{split}
\end{equation}
\end{Beweis}

We further provide the following specializations of the ubiquitous
function $N\in\R^\X\!$,
$N(x)=\prod_{j\in(n]}N_j(x_j)$\´:\pagebreak[2]

\begin{Lemma}\label{eqN}
  Let $E_l:=\{\,c\in\R\mit c^{l}=1\,\}$ denote the set
  of the \(l^\textit{th}\) roots of unity in $\R$.
  For $x\,\in\,\X_j\sb\R$ hold\´:
  \begin{enumerate}[(i)]
    \item If $\X_j=E_{d_j+1}$ (\´$\abs{E_{d_j+1}}=d_j+1$\´) and\\
     if $\R$ is an integral domain\´:\ \hfill
     \framebox[13em][l]{\strut$N_j(x)\,=\,(d_j+1)\,x^{-1}$.}
    \item If $\X_j\uplus\{0\}$ is a finite subfield of $\R$\´:\ \hfill
     \framebox[13em][l]{\strut$N_j(x)\,=\,-x^{-1}$.}
    \item \parbox[t]{21em}{\strut If $\X_j=E_{d_j}\uplus\{0\}$ (\´$\abs{E_{d_j}}=d_j$\´)
     \smallskip and\\
     if $\R$ is an integral domain\´:}\\[-2.1\baselineskip]\null\ \hfill
     \framebox[13em][l]{$N_j(x)\,=\,\begin{cases}
        d_j1 & \text{for}\ x\neq0\,,\\
        -1 & \text{for}\ x=0\,.
      \end{cases}$}
    \item If $\X_j$ is a finite subfield
      of $\R$\´:\ \hfill
     \framebox[13em][l]{\strut$N_j(x)\,=\,-1$ .}
    \item If $\X_j=\{0,1,\dotsc,d_j\}\sb\Z$\´:\ \hfill
     \framebox[13em][l]{\strut$N_j(x)
      \,=\,(-1)^{d_j+x}\,d_j!\,\binom{d_j}{x}^{\!\!-1}\!$.}
    \item For $\alpha\in\R$ we have\´:\ \hfill
     \framebox[13em][l]{\strut$N_{\X_j+\alpha}(x+\alpha)\,=\,N_{\X_j}(x)$
     .}
  \end{enumerate}
\end{Lemma}\vsp

\begin{Beweis}
For finite subsets $\D\sb\R$ we define
\begin{equation}\ps
L_\D(X)\,:=\,\prod_{\hx\in\D}(X-\hx)\,\ .
\end{equation}
It is well"=known that, if $E_l$ contains $l$ elements and lies in
an integral domain,
\begin{equation}\ps
L_{E_l}(X)
  \,=\,\prod_{\hx\in E_l}(X-\hx)
  \,=\,X^{l}-1
  \,=\,(X-1)(X^{l-1}+\dotsb+X^0)\,\ .
\end{equation}
Thus
\begin{equation}\ps
  L_{E_l\!\ssm1}(1)
  \,=\,\frac{\prod_{\hx\in E_l}(X-\hx)}{X-1}\bigm|_{X=1}
  \,=\,\frac{X^{l}-1}{X-1}\bigm|_{X=1}
  \,=\,X^{l-1}+\dotsb+X^0\bigm|_{X=1}
  \,=\,l1\ .
\end{equation}
Using this, 
we get for $x\in E_l$
\begin{equation}\ps\label{NEl}
  L_{E_l\!\ssm x}(x)
  \,=\,L_{x(E_l\!\ssm1)}(x)
  \,=\prod_{\hx\in E_l\ssm1\!\!\!}(x-x\hx)
  \,=\,x^{l-1}L_{E_l\!\ssm1}(1)
  \,=\,l x^{-1}\ .
\end{equation}
This gives \((i)\) with $l=\abs{\X_j}=d_j+1$.

Part\,\((ii)\) is a special case of part\,\((i)\), where
$\X_j\,=\,F_{p^k}\!\ssm0\,=\,E_{p^k-1}$ and where consequently
$d_j+1=\abs{\X_j}=(p^k-1)\equiv-1\pmod{p}$.

To get $N_j(x)=L_{\{0\}\uplus E_l\!\ssm x}(x)$ with $x\neq0$ in
part\,\((iii)\) and part\,\((iv)\) we multiply Equation\,\eqref{NEl}
with $x-0$ and use $l=\abs{\X_j}-1=p^k-1\equiv-1\pmod{p}$ for
part\,\((iv)\) and $l=\abs{\X_j}-1=d_j$ for part\,\((iii)\).
For $x=0$ 
we obtain in part\,\((iii)\) and part\,\((iv)\)
\begin{equation}\ps
N_j(0)
 \,=\,L_{E_l}(0)
 \,=\,\prod_{\hx\in E_l}(-\hx)
 \,=\,-\!\!\!\!\!\!\prod_{\ \ \hx\in E_l\!\sm\{1,-1\}\!\!\!\!\!\!\!\!}\!\!\!\!(-\hx)
 \,=\,-1\ \ ,
\end{equation}
since each subset $\{\hx,\hx^{-1}\}\sb E_l\!\sm\{1,-1\}$ contributes
$(-\hx)\,(-\hx^{-1})=1$ to the product \hbox{-- \,as}
$\hx\neq\hx^{-1}\!$, since $\hx^2-1=0$ holds only for
$\hx=\pm1\!{}\´$ -- \´\,and $E_l\!\sm\{1,-1\}$ is partitioned by
such subsets. This completes the proofs of parts\,\((iii)\) and
\((iv)\).

We now turn to part\,\((v)\)\´:
\begin{equation}\ps\label{N012}
 N_j(x)
 \,=\,\bigl(\!\!\prod_{0\leq\hx<x}\!\!(x-\hx)\,\bigr)\!\prod_{x<\hx\leq d_j}\!\!(x-\hx)
 \,=\,x!\,(d_j-x)!\,(-1)^{d_j-x}
 \,=\,(-1)^{d_j+x}\,d_j!\,\binom{d_j}{x}^{\!\!-1}\!.
\end{equation}

Part\,\((vi)\) is trivial.
\end{Beweis}\vspp


\section{Interpolation polynomials and inversion formulas}\label{sec.ip}

\index{$\p$, $\R[X^{\leq d}]$, $\R^{[d]}$, $\R^\X$} This section may
be skipped at a first reading; the only things you need from here to
understand the rest of the paper are\´:\smallskip\\
--\hbox{\ \ }the fact that grids
$\X:=\X_1\times\dotsb\times\X_n\,\sb\,\R^n$ over integral domains
      $R$ are always
\phantom{--\hbox{\ \ }}\emph{integral grids},
      in the sense of Definition\,\ref{sz.dvg}\´, and\\
--\hbox{\ \ }the inversion formula\,\´\ref{sz.itg}\,, which is,
      in this case, just the well"=known interpolation
\phantom{--\hbox{\ \ }}formula for polynomials
      applied to polynomial maps $P|_{\X}$.\smallskip\\
The rest of this section is concerned with providing some generality that is only needed in
the very last section of this paper.\bigskip

We have to investigate the canonical homomorphism \rand$"\p"\DP
P\lmto P|_\X$ that maps polynomials $P$ to polynomial maps
 $P|_\X\DP x\mto P(x)$ on a fixed \(d\)"~grid $\X\sb\R^n\!$.
As the monic polynomial
\rand$"L_j"=L_{\X_j}(X_j):=\prod_{\hx\in\X_j}(X_j-\hx)$ maps all
elements of $\X_j$ to $0$, we may replace each given polynomial $P$
by any other polynomial of the form $P+\sum_{j\in(n]}H_jL_j$ without
changing its image $P|_\X$. By applying such modifications, we may
assume that $P$ has partial degrees $\deg_j(P)\leq \abs{\X_j}-1=d_j$
(\´see Example\,\ref{bsp.cn} for an illustration of this method\´).
Hence the image of $\p$ does not change if we regard $\p$ as a map
on $\R[X^{\leq d}]$ (\´which we identify with $\R^{[d]}$ by
$P\mto(P_\d)_{\d\in[d]}$\´). The resulting map
 \rand\begin{equation}
"\p"\DP\R[X^{\leq d}]=\R^{[d]}\lto\R^\X\,,\,\
       P\lmto P|_\X\,:=\,(\´x\mto P(x)\´)
\end{equation}
is in the most important cases an isomorphism or at least a
monomorphism, as we will see in this section. In general, however,
the situation is much more complicated, we give a short example and
make a related, more general remark\´:

\begin{Beispiel}
Over $\R=\Z_6:=\Z/6\Z$ we have $X^3|_{\Z_6}=X|_{\Z_6}$ and
$3X^2|_{\Z_6}=3X|_{\Z_6}$, so that each polynomial map
$\X:=\Z_6\lto\Z_6$ can be represented by a polynomial of the form
$aX^2+bX+c$, with $a\in\{0,1,-1\}$. Hence the corresponding
$3\cdot6^2$ distinct maps are the only maps out of the $6^6$ maps
from $\X=\Z_6$ to $\Z_6$ that can be represented by polynomials at
all. This simple example shows also that the kernel $\ker(\p)$ may
be very complicated even in just one dimension.
\end{Beispiel}

\begin{Bemerkung}
There are some general results for the rings $\R=\Z_m$
of integers mod $\!m$:\smallskip\\
--\hbox{\ \ }In \cite{must} a system of polynomials in
  $\Z_m[X_1,\dots,X_n]$ is given that represent all poly-
\phantom{--\hbox{\ \ }}nomial maps
$\Z_m\!\n\!\lto\Z_m$ and the number of all such maps is determined.\\[2pt]
--\hbox{\ \ }In \cite{sp} it is shown that the Newton algorithm can
  be used to determine interpolation
\phantom{--\hbox{\ \ }}polynomials, if they exist. The ``divided
  differences'' in this algorithm are, like the
\phantom{--\hbox{\ \ }}interpolation polynomials themselves, not
  uniquely determined over arbitrary commu-
\phantom{--\hbox{\ \ }}tative rings, and exist if and only if
interpolation polynomials exist.
\end{Bemerkung}\bigskip

But back to the main subject. In which situations does $\p\DP P\lmto
P|_\X$ become an isomorphism, or equivalently, when does its
representing matrix $\Phi$ possess an inverse? Over commutative
rings $\R$, square matrices $\Phi\in\R^{m\times m}$ with
nonvanishing determinant do not have an inverse, in general.
However, there is the matrix \rand$"\Adj(\Phi)"$ -- the adjoint or
cofactor matrix -- that comes close to being an inverse\´:
\begin{equation}
\Phi\Adj(\Phi)\,=\,\Adj(\Phi)\´\Phi\,=\,\det(\Phi)\´\E\,\ .
\end{equation}
In our concrete situation, where \rand$"\Phi"\in\R^{\X\times[d]}$ is
the matrix of $\p$ (\´a tensor product of Vandermonde matrices\´),
we work with \rand$"\Psi"$ (\´from Definition\,\ref{def.N}\´)
instead of the adjoint matrix $\Adj(\Phi)$. $\Psi$ comes closer
than $\Adj(\Phi)$ to being a right inverse of $\Phi$. %
The following theorem shows that
\begin{equation}
\Phi\Psi\,=\,\bigl(N(x)\,\?{\tx=x}\bigr)_{\tx,x\in\X}\,\ ,
\end{equation}
and the entries $N(x)$ of this diagonal matrix divide the entries
$\det(\Phi)$ of $\Phi\Adj(\Phi)$, so that $\Phi\Psi$ is actually
closer than $\Phi\Adj(\Phi)$ to the unity matrix (\´provided we
identify the column indices $x\in\X$ and row indices $\d\in[d]$ in
some way with the numbers $1,2,\dotsc,\abs{\X}=\abs{[d]}$, in order
to make $\det(\Phi)$ and $\Adj(\Phi)$ defined\´).

However, we used the matrix
 \randd$"\Phi"\in\R^{\X\times[d]}$ of \randd$"\p"\DP P\lmto P|_\X$
here just to explain the role of $\Psi$. In what follows, we do not
use it any more; rather, we prefer notations with ``$\p$'' or
``$|_\X$.'' For maps/tuples $y\in\R^\X\!$, we write
\rand$"(\Psi y)(X)"\in\R[X^{\leq d}]$, as already defined, for the
polynomial whose coefficients form the tuple
$\Psi y\in\R^{[d]}\!$, i.e., $(\Psi y)(X)=\Psi y$ by identification.
We have\´:

\index{interpolation}
\begin{Satz}[Interpolation]\label{sz.ip}
For maps $y\DP\X\lto\R$,
$$\boxed{\strut\,(\Psi y)(X)|_\X\,=\,Ny\,}\ .$$
\end{Satz}\vspp

\begin{Beweis}
As both sides of the equation are linear in $y$, it suffices to
prove the equation for the maps $y=e_\tx$, where $\tx$ ranges over
$\X$. Now we see that, at each point $x\in\X$, we actually have
\begin{equation}\ps
 (\Psi e_\tx)(X)|_\X(x)
 \ \eqby{\ref{lem.lagr}}\ L_{\X,\tx}(x)
 \ =\ N(x)\,\?{x=\tx}
 \ =\ (Ne_\tx)(x)\ .\vspace{-.5\baselineskip}
\end{equation}
\end{Beweis}

With this theorem, we are able to characterize the situations in
which $\p\DP P\lmto P|_\X$ is an isomorphism\´:

\index{grid!division}
\begin{EquiDef}[Division grids]\label{eqdef.dvg}
  We call a \(d\)"~grid $\X\sb\R^n$ a \emph{division grid} (\´over $\R$\´)
  if it has the following equivalent properties\´:\pagebreak[2]
\begin{enumerate}[(i)]
\item For all $j\in(n]$ and all $x,\tx\in\X_j$ with $x\neq\tx$
      the difference $x-\tx$ is invertible.
\item $N=N_\X$ is pointwise invertible, i.e., for all $x\in\X\!$, $N(x)$ is invertible.
\item $\Pi N$ is invertible.
\item $\p\DP\R[X^{\leq d}]=\R^{[d]}\lto\R^\X$\! is bijective.
\end{enumerate}
\end{EquiDef}\vsp

\begin{Beweis}
The equivalence of \((i)\),\((ii)\) and \((iii)\) follows from the
Definition\,\ref{def.N} of $N$\!, the definition
 $\Pi N=\prod_{x\in\X}N(x)$ and the associativity and commutativity
of $\R$.

Assuming \((ii)\), it follows from Theorem\,\ref{sz.ip} that
$y\lmto(\Psi(N^{-1}y))(X)$ is a right inverse of
 $\p\DP P\lmto P|_{\X}$\,:
\begin{equation}\ps
y\,\lmto\,(\Psi(N^{-1}y))(X)\eqby[\,\lmto\,]{\p}N(N^{-1}y)=y\,\ .
\end{equation}
It is even a two"=sided inverse, since square matrices $\Phi$ over a
commutative ring $\R$ are invertible from both sides if they are
invertible at all (\´since $\Phi\Adj(\Phi)=\det(\Phi)\´\E$\´). This
gives \((iv)\).

Now assume \((iv)\) holds; then for all $x\in\X$,
\begin{equation}\ps
 \bigl(\psi_{\d,x}\bigr)_{\d\in[d]}
 \ =\ \Psi e_x
 \ \eqby{\ref{sz.ip}}\ \p^{-1}(Ne_x)
 \ =\ N(x)\,\p^{-1}(e_x)\,\ ,
\end{equation}
and in particular,
\begin{equation}\ps
 1\,\eqby{\eqref{eqP2}}\,\psi_{d,x}
  \,=\,N(x)\,\bigl(\p^{-1}(e_x)\bigr)_{\!d}\,\ .
\end{equation}
Thus the $N(x)$ are invertible and that is \((ii)\).
\end{Beweis}

If $\p\DP\R[X^{\leq d}]\lto\R^\X$ is an isomorphism, then
 \rand$"\p^{-1}"(y)$
is the unique polynomial in $\R[X^{\leq d}]$ that interpolates a
given map $y\in\R^\X\!$, so that, by Theorem\,\ref{sz.ip}\,, it has
to be the polynomial $\Psi(N^{-1}y)\in\R^{[d]}=\R[X^{\leq d}]$. This
yields the following result\´:

\index{interpolation formula}\index{grid!Boolean}
\begin{Satz}[Interpolation formula]\label{sz.dvg}
Let $\X$ be a division grid (\´e.g., if $\R$ is a field or if $\X$
is the \emph{Boolean grid} $\{0,1\}^n$). For $y\in\R^\X\!$,
$$\boxed{\strut\,\p^{-1}(y)\,=\,\Psi(N^{-1}y)\,}\ .$$
\end{Satz}\vspp

\noindent This theorem can be found in \cite[Theorem\,2.5.2]{da},
but just for fields $\R$ and in a different representation (\´with
$\p^{-1}(y)$ as a determinant\´).\pagebreak[2]

Additionally, if $\X$ is not a division grid, we may apply the
canonical localization homomorphism
 \index{$\pi$, $\S$, $\R_N$}\Rand{\(\pi\)}\Rand{\(\S\´,\,\R_N\)}
\begin{equation}
 \pi\DP\R\lto\R_N:=\S^{-1}\R\ ,\,\
         r\lmto r^\pi :=\tfrac{r}{1}
        \ \quad\text{with}\quad
 \S:=\{\,(\Pi N)^m\mit m\in\N\,\}\,\ ,
\end{equation}
and exert our theorems in this situation. As $\pi$ and $\R_N$ have
the universal property with respect to the invertibility of $(\Pi
N)^\pi$ in $\R_N$ (\´as required in \ref{eqdef.dvg}\((iii)\)\´),
$\pi$ and $\R_N$ are the best choices. This means specifically that
if $(\Pi N)^\pi$ is not invertible in the codomain $\R_N$ of $\pi$,
then no other homomorphism $\pi'$ has this property, either. In
general, $\pi$ does not have this property itself\´: By definition,
\begin{equation}
\frac{r_1}{s_1}\,=\,\frac{r_2}{s_2}\ \ :\lEqi\ \ \ex s\in\S\DP\
s\,r_1s_2\,=\,s\,r_2s_1\,\ ,
\end{equation}
 and hence
\rand\begin{equation}
"\ker(\pi)"
           \,=\,\{\,r\in\R\mit\ex m\in\N\DP(\Pi N)^m\,r=0\,\}\,\ ,
\end{equation}
 so that $(\Pi N)^\pi=0$ is possible. Localization works in the
following situation\´:

\index{grid!affine}
\begin{EquiDef}[Affine grids]\label{eqdef.afg}
  We call a \(d\)"~grid $\X\sb\R^n$ \emph{affine} (\´over $\R$\´)
  if it has the following equivalent properties\´:
\begin{enumerate}[(i)]
\item $\Pi N$ is not nilpotent.
\item $\pi\neq0$.
\item $(\Pi N)^\pi$ is invertible in $\R_N$.
\item $\pi\neq0$ is injective on the $\X_j$\\ and hence induces a bijection
  \rand$\X\lto"\X^\pi":={\X_1}^{\!\!\pi}\times\dotsb\times
  {\X_n}^{\!\!\!\pi}$.
\end{enumerate}
\end{EquiDef}\vsp

\begin{Beweis}
Part\,\((ii)\) is equivalent to $1^\pi\neq0$, and this means that
$s1\neq0$ for all $s$ in the multiplicative system $\S=\{\,(\Pi
N)^m\mit m\in\N\,\}$; thus $(i)\lEqi(ii)$.

Of cause \((\Pi N)^\pi\frac{1}{\Pi N}=\frac{1}{1}\ms\) is the unity
in $\R_N$, provided $\frac{1}{1}=1^\pi\neq0$; thus $(ii)\lTo(iii)$.

If \((iii)\) holds then $(\Pi N)^\pi$ and its factors
$(x_j-\tx_j)^\pi$ do not vanish; thus $(iii)\lTo(iv)$.

Finally, the implication $(iv)\lTo(ii)$ is trivial.
\end{Beweis}

If $\X\sb\R^n$ is affine, then
 \rand$"\X^\pi":={\X_1}^{\!\!\pi}\times\dotsb\times {\X_n}^{\!\!\!\pi}\sb \R_N\n[\!\!n]$
is a division \(d\)"~grid over $\R_N$ by \ref{eqdef.afg}\,\((iv)\),
\ref{eqdef.afg}\,\((iii)\) and \ref{eqdef.dvg}\,\((iii)\). Now,
Theorem\,\ref{sz.dvg} applied to \rand$y:="P^\pi"|_{\X^\pi}$ with
$P^\pi=\sum_{\d\in[d]}{P_\d}^{\!\pi} X^\d$ yields
\begin{equation}
P^\pi\,=\,\Psi_{\X^\pi}\bigl((N_{\X^\pi})^{-1}(P^\pi|_{\X^\pi})\bigr)\,\ ,
\end{equation}
along with the associated constants
$N_{\X^\pi}\in\R_N\n[\!\!\X^\pi]\!$ and
$\Psi_{\X^\pi}\in\R_N^{\,[d]\times\X^\pi}\!\!$ of
$\X^\pi\!$.\pagebreak[2]

With componentwise application of $\pi\DP r\mto \tfrac{r}{1}$ to
$P|_\X,\,N\in\R^\X$ and to $\Psi\in\R^{[d]\times\X}$ so that
\randd$(P|_\X)^\pi\!,\,"N^\pi"\in\´\R_N\n[\!\!\X]$ and
\randd$"\Psi^\pi"\in\´\R_N^{\,[d]\times\X}\!$, we obtain\´:

\index{inversion formula}
\begin{Satz}[Inversion formula]\label{sz.afg}
Let $\X$ be affine (\´e.g., if $\R$ does not possess nilpotent
elements\´). For $P\,\in\,\R[X^{\leq d}]=\R^{[d]}\!$,
$$\boxed{\strut\,P^\pi\,=\,\Psi^\pi\bigl((N^\pi)^{-1}(P|_\X)^\pi\bigr)\,}\ .$$
\end{Satz}\vspp

If $\pi$ is injective on its whole domain $\R$ then $\R$ is a
subring of $\R_N$ and we may omit $\pi$ in formula\,\ref{sz.afg}\,.
In fact, we will see that this is precisely when $\p$ is
injective, as seen in the following characterization\´:

\index{grid!integral}
\begin{EquiDef}[Integral grids]\label{eqdef.itg}
  We call a \(d\)"~grid $\X\sb\R^n$ \emph{integral} (\´over $\R$\´)
  if it has the following, equivalent properties\´:
\begin{enumerate}[(i)]
\item For all $j\in(n]$ and all $x,\tx\in\X_j$
  with $x\neq\tx\!$, $x-\tx$ is not a zero divisor.
\item For all $x\in\X\!$, $N(x)$ is not a zero divisor.
\item $\Pi N$ is not a zero divisor.
\item $\pi$
  is injective ($\R\sb\R_N$).
\item $\p\DP\R[X^{\leq d}]=\R^{[d]}\lto\R^\X$ is injective.
\end{enumerate}
\end{EquiDef}\vsp

\begin{Beweis}
The equivalence of \((i)\),\((ii)\) and \((iii)\) follows from the
Definition\,\ref{def.N} of $N$, the definition
 $\Pi N=\prod_{x\in\X}N(x)$ and the associativity and  commutativity
of $\R$.

  As already mentioned $\ker(\pi)=\{\,r\in\R\mit\ex m\in\N\DP(\Pi N)^mr=0\,\}$,
  so $(iii)\lTo(iv)$.

  If \((iv)\) holds, then $\Pi N$ is invertible in $\R_N$. By
  Equivalence\,\ref{eqdef.dvg}\,,
  it follows that $\p\DP\R_N[X^{\leq d}]\lto\R_N\n[\!\!\X]$ is bijective,
  so that $(iv)\lTo(v)$.

  Now suppose that \((ii)\) does not hold, so that there are a point $x\in\X$
  and a constant $M\in\R\ssm0$ with
\begin{equation}\ps
MN(x)\,=\,0\ .
\end{equation}
Then
\begin{equation}\ps
P\,:=\,\Psi(Me_x)\,\neq\,0\,\ ,
\end{equation}
as
\begin{equation}\ps
P_d\,=\,M(\Psi(e_x))_d\,=\,M\psi_{d,x}\,\eqby{\eqref{eqP2}}M\,\neq\,0\,\ .
\end{equation}
However,
\begin{equation}\ps
\p(P)\,\eqby{\ref{sz.ip}}\,NMe_x\,=\,MN(x)e_x\,\equiv\,0\,\ ,
\end{equation}
so that \((v)\) does not hold, either. Thus $(v)\lTo(ii)$.
\end{Beweis}\pagebreak[2]

Any integral grid $\X$ over $\R$ is, in fact, a division grid over
$\R_N\sp\R$, since $\Pi N$ becomes invertible in $\R_N$.
Formula\,\ref{sz.dvg} applied to $y:=P|_\X$ yields the following
specialization of Theorem\,\ref{sz.afg}\´:\pagebreak[2]

\index{inversion formula}
\begin{Satz}[Inversion formula]\label{sz.itg}
Let $\X$ be integral (\´e.g., if $\R$ is an integral domain\´). For
$P\,\in\,\R[X^{\leq d}]=\R^{[d]}\!$,
$$\boxed{\strut\,P\,=\,\Psi(N^{-1}P|_\X)\,}\ .$$
\end{Satz}\vspp

From the case $P=1$, we see that $N^{-1}P|_\X$ inside this formula
does not lie in $\R^\X$ in general (\´of course
$N^{-1}P|_\X\in\R_N\n[\!\!\X]$\´). This also shows that, in general,
the maps of the form $Ny$, with $y\in\R^\X\!$, in
Theorem\,\ref{sz.ip} are not the only maps that can be represented
by polynomials over $\R$, i.e.,
 $\{\,Ny\mit y\in\R^\X\,\}\varsubsetneq\Bild(\p)$.
However, the maps of the form $Ny$ are exactly the linear
combinations of Lagrange's polynomial maps $Ne_x=L_{\X,x}|_\X$ over
the grid $\X$; and if we view, a bit more generally, Lagrange
polynomials $L_{\tilde\X,x}$ over subgrids
$\tilde\X=\tilde\X_1\times\dotsb\times\tilde\X_n\sb\X$, then the
maps of the form $L_{\tilde\X,x}|_\X$ span $\Bild(\p)$, as one can
easily show.

On the other hand, in general, $\Bild(\p)\varsubsetneq\R^\X\!$, so
that not every map $y\in\R^\X$ can be interpolated over $\R$. If
$\X$ is integral, then interpolation polynomials exist over the
bigger ring $\R_N$. The univariate polynomials
$\tbinom{X}{k}:=\tfrac{X(X-1)\dotsb(X-k+1)}{k!}$, for example,
describe integer"=valued maps (\´on the whole domain $\Z$\´), but do
not lie in $\Z[X]$. More information about such ``overall''
integer"=valued polynomials over quotient fields can be found, for
example, in \cite{ccf} and \cite{ccs}, and in the literature cited
there.\medskip

\pagebreak[2]The reader might find it interesting that the principle
of inclusion and exclusion follows from Theorem\,\ref{sz.itg} as a
special case\´:

\index{inclusion and exclusion}
\begin{Proposition}[Principle of inclusion and exclusion]\label{prop.ip}\*\\[1pt]
Let $\X:=\{0,1\}^n=[d]$ and $x\in\X$; then $x^\d=\?{\d\leq x}$ for
all $\d\in[d]$. Thus, for arbitrary
 $P=(P_\d)
 \,\in\,\R^{[d]}=\R[X^{\leq d}]$,
\begin{equation}\label{eq.M}
  P(x)\ =\ \smash{\sum\nolimits_{\d\leq x}}\,P_\d\,\ .\vsp
\end{equation}
Formula\,\ref{sz.itg} is the Möbius inversion to
Equation\,\eqref{eq.M}\´:
\begin{equation}
\begin{split}
  P_\d
  \ceqby{\ref{sz.itg}}\sum_{x\in[d]}\,\psi_{\d,x}N^{-1}(x)P(x)\\
  \ceqby{\ref{def.N}}\sum_{x\in[d]}\,
     \bigl[{\textstyle\prod\limits_{j\in(n]}}\´\?{x_j\leq\d_j}\,(-1)^{1-\d_j}\´\bigr]
     \bigl[{\textstyle\prod\limits_{j\in(n]}}\´(-1)^{1-x_j}\´\bigr]P(x)\\
  \ceq\sum\nolimits_{x\leq\d}\,(-1)^{\Sigma(\d-x)}P(x)\,\ .
\end{split}
\end{equation}
\end{Proposition}\vspp\clearpage


\section{Coefficient formulas -- the main results}\label{sec.cf}

The applications in this paper do not start with a map $y\in\R^\X$
that has to be interpolated by a polynomial $P$. Rather, we start
with a polynomial $P\!$, or with some information about a polynomial
$P\in\R[X]$, which describes the very map $y:=P|_\X$ that we would
like to understand. Normally, we will not have complete information
about $P\!$, so that we do not usually know all coefficients $P_\d$
of $P$. However, there may be a coefficient $P_\d$ in
$P=\sum_{\d\in\N^n}P_\d X^\d$ that, on its own, allows conclusions
about the map $P|_\X$. We define (\´see also figure\,1 below\´)\´:

\index{\(d\)-leading}
\begin{Definition}\label{def.dl}
Let $P=\sum_{\d\in\N^n}P_\d X^\d\in\R[X]$ be a polynomial. We call a
multiindex $\e\leq d\in\N^n$ 
\emph{\(d\)"=leading} in $P$ if for each monomial $X^\d$ in $P\!$,
i.e., each $\d$ with $P_\d\neq0$, holds either\smallskip

--\ $\displaystyle
 \text{(case\,1) \,$\d=\e$;\quad \,\emph{or}}$

--\ $\displaystyle
 \text{(case\,2) \,\,there is a ${j\in(n]}$ such that
 $\d_j\neq\e_j$ but $\d_j\leq d_j$\,\ .}$
\end{Definition}\vspp

Note that the multiindex $d$ is \(d\)"=leading in polynomials $P$
with $\deg(P)\leq\Sigma d$. In this situation, case\,2 reduces to
``there is a $j\in(n]$ such that $\d_j<d_j$,\!'' and, as
$\Sigma\d\leq\Sigma d$ for all $X^\d$ in $P\!$, we can conclude\´:
\begin{equation}
\text{``not case\,2''}
 \,\ \lTo\ \ \d\geq d
 \ \ \lTo\ \ \d=d
 \ \ \lTo\,\ \text{``case\,1''}\ .
\end{equation}
Thus $d$ really is \(d\)"~leading in $P$ (\´see also figure\,2 on
page\,\pageref{Fig2}\´). Of course, if all partial degrees are
restricted by $\deg_j(P)\leq d_j$ then all multiindices $\d\leq d$
are \(d\)"=leading. Figure\,1 (\´below\´) shows a nontrivial example
$P\in\R[X_1,X_2]$. The monomials $X^\d$ of $P$ ($P_\d\neq0$), and
the $2^n\!-1=3$ ``forbidden areas'' of each of the two
\(d\)"=leading multiindices, are marked.
\afterpage{%
\begin{figure}
\begin{center}
\setlength{\unitlength}{0.00087489in}
\begingroup\makeatletter\ifx\SetFigFont\undefined
\def\x#1#2#3#4#5#6#7\relax{\def\x{#1#2#3#4#5#6}}%
\expandafter\x\fmtname xxxxxx\relax \def\y{splain}%
\ifx\x\y   
\gdef\SetFigFont#1#2#3{%
  \ifnum #1<17\tiny\else \ifnum #1<20\small\else
  \ifnum #1<24\normalsize\else \ifnum #1<29\large\else
  \ifnum #1<34\Large\else \ifnum #1<41\LARGE\else
     \huge\fi\fi\fi\fi\fi\fi
  \csname #3\endcsname}%
\else
\gdef\SetFigFont#1#2#3{\begingroup
  \count@#1\relax \ifnum 25<\count@\count@25\fi
  \def\x{\endgroup\@setsize\SetFigFont{#2pt}}%
  \expandafter\x
    \csname \romannumeral\the\count@ pt\expandafter\endcsname
    \csname @\romannumeral\the\count@ pt\endcsname
  \csname #3\endcsname}%
\fi
\fi\endgroup
{\renewcommand{\dashlinestretch}{30}
\begin{picture}(3400,2500)(0,-350)
\put(2270,390){\shade\ellipse{120}{120}}
\put(1970,390){\shade\ellipse{120}{120}}
\put(1970,990){\shade\ellipse{120}{120}}
\put(1670,390){\blacken\ellipse{120}{120}}
\put(1670,690){\blacken\ellipse{120}{120}}
\put(1670,990){\blacken\ellipse{120}{120}}
\put(1670,1590){\shade\ellipse{120}{120}}
\put(1370,690){\blacken\ellipse{120}{120}}
\put(1370,990){\blacken\ellipse{120}{120}}
\put(1370,1290){\shade\ellipse{120}{120}}
\put(1070,390){\blacken\ellipse{120}{120}}
\put(1070,990){\blacken\ellipse{120}{120}}
\put(770,390){\blacken\ellipse{120}{120}}
\put(770,690){\blacken\ellipse{120}{120}}
\put(770,1290){\shade\ellipse{120}{120}}
\put(770,1590){\shade\ellipse{120}{120}}
\put(470,390){\blacken\ellipse{120}{120}}
\put(470,990){\blacken\ellipse{120}{120}}
\path(320,390)(320,2040)
\path(368,1880)(320,2040)(272,1880)
\path(320,1590)(305,1590)
\path(320,1290)(305,1290)
\path(320,990)(305,990)
\path(320,690)(305,690)
\path(320,390)(305,390)
\path(470,240)(2720,240)
\path(2560,192)(2720,240)(2560,288)
\path(2270,240)(2270,225)
\path(1970,240)(1970,225)
\path(1670,240)(1670,225)
\path(1370,240)(1370,225)
\path(1070,240)(1070,225)
\path(770,240)(770,225)
\path(470,240)(470,225)
\thicklines
\put(1070,690){\ellipse{120}{120}}
\put(470,690){\ellipse{120}{120}}
\allinethickness{.6pt}
\put(995,1215){\arc{90}{1.5708}{3.1416}}
\put(1145,1215){\arc{90}{0}{1.5708}}
\path(950,1215)(950,1765)
\path(1190,1765)(1190,1215)
\path(1145,1170)(995,1170)
\put(395,1215){\arc{90}{1.5708}{3.1416}}
\put(545,1215){\arc{90}{0}{1.5708}}
\path(350,1215)(350,1765)
\path(590,1765)(590,1215)
\path(545,1170)(395,1170)
\put(1895,615){\arc{90}{1.5708}{3.1416}}
\put(1895,765){\arc{90}{3.1416}{4.7124}}
\path(1850,615)(1850,765)
\path(1895,810)(2445,810)
\path(2445,570)(1895,570)
\put(1895,1215){\arc{90}{1.5708}{3.1416}}
\path(1850,1215)(1850,1765)
\path(2445,1170)(1895,1170)
\thinlines
\dashline{96.000}(1820,240)(1820,1140)
\dashline{96.000}(320,1140)(1790,1140)
\path(550,1170)(590,1210)
\path(450,1170)(590,1310)
\path(365,1185)(590,1410)
\path(350,1270)(590,1510)
\path(350,1370)(590,1610)
\path(350,1470)(590,1710)
\path(350,1570)(545,1765)
\path(350,1670)(445,1765)
\path(1150,1170)(1190,1210)
\path(1050,1170)(1190,1310)
\path(965,1185)(1190,1410)
\path(950,1270)(1190,1510)
\path(950,1370)(1190,1610)
\path(950,1470)(1190,1710)
\path(950,1570)(1145,1765)
\path(950,1670)(1045,1765)
\path(2350,570)(2445,665)
\path(2250,570)(2445,765)
\path(2150,570)(2390,810)
\path(2050,570)(2290,810)
\path(1950,570)(2190,810)
\path(1865,585)(2090,810)
\path(1850,670)(1990,810)
\path(1850,770)(1890,810)
\path(2350,1170)(2445,1265)
\path(2250,1170)(2445,1365)
\path(2150,1170)(2445,1465)
\path(2050,1170)(2445,1565)
\path(1950,1170)(2445,1665)
\path(1865,1185)(2445,1765)
\path(1850,1270)(2345,1765)
\path(1850,1370)(2245,1765)
\path(1850,1470)(2145,1765)
\path(1850,1570)(2045,1765)
\path(1850,1670)(1945,1765)
\put(0,-150){\makebox(0,0)[lb]%
{\smash{{{\SetFigFont{10}{12.0}{sl}Figure\,1\´: Monomials of a polynomial $P$ with}}}}}
\put(0,-350){\makebox(0,0)[lb]%
{\smash{{{\SetFigFont{10}{12.0}{sl}\((4,2)\)-leading multiindices $(0,1)$ and $(2,1)$.}}}}}
\put(-20,1820){\makebox(0,0)[lb]{\smash{{{\SetFigFont{8}{9.6}{rm}$\deg_2$}}}}}
\put(220,1560){\makebox(0,0)[lb]{\smash{{{\SetFigFont{8}{9.6}{rm}4}}}}}
\put(220,1260){\makebox(0,0)[lb]{\smash{{{\SetFigFont{8}{9.6}{rm}3}}}}}
\put(220,960){\makebox(0,0)[lb]{\smash{{{\SetFigFont{8}{9.6}{rm}2}}}}}
\put(220,660){\makebox(0,0)[lb]{\smash{{{\SetFigFont{8}{9.6}{rm}1}}}}}
\put(220,360){\makebox(0,0)[lb]{\smash{{{\SetFigFont{8}{9.6}{rm}0}}}}}
\put(2420,80){\makebox(0,0)[lb]{\smash{{{\SetFigFont{8}{9.6}{rm}$\deg_1$}}}}}
\put(2235,120){\makebox(0,0)[lb]{\smash{{{\SetFigFont{8}{9.6}{rm}6}}}}}
\put(1935,120){\makebox(0,0)[lb]{\smash{{{\SetFigFont{8}{9.6}{rm}5}}}}}
\put(1635,120){\makebox(0,0)[lb]{\smash{{{\SetFigFont{8}{9.6}{rm}4}}}}}
\put(1335,120){\makebox(0,0)[lb]{\smash{{{\SetFigFont{8}{9.6}{rm}3}}}}}
\put(1035,120){\makebox(0,0)[lb]{\smash{{{\SetFigFont{8}{9.6}{rm}2}}}}}
\put(735,120){\makebox(0,0)[lb]{\smash{{{\SetFigFont{8}{9.6}{rm}1}}}}}
\put(435,120){\makebox(0,0)[lb]{\smash{{{\SetFigFont{8}{9.6}{rm}0}}}}}
\end{picture}
}
\end{center}
\end{figure}}

In what follows, we examine how the preconditions of the inversion
formula\,\ref{sz.itg} may be weakened. It turns out that
formula\,\ref{sz.itg} holds without further degree restrictions for
the \(d\)"~leading coefficients $P_\e$ of $P$\!. The following
theorem is a generalization and a sharpening of Alon and Tarsi's
(\´second\´) Combinatorial Nullstellensatz
\cite[Theorem\,1.2]{al2}\´:

\index{coefficient formula}
\begin{Satz}[Coefficient formula]\label{sz.cn}
  Let $\X
$ be an integral \(d\)"~grid.
  For each polynomial $P=\sum_{\d\in\N^n}P_\d X^\d\,\in\,\R[X]$
  with \(d\)"=leading multiindex $\e\leq d\in\N^n$\!,
\begin{enumerate}[(i)]
\item \framebox[10em][l]{\strut$P_\e\,=\,(\Psi(N^{-1}P|_\X)\´)_\e$}\quad
  (\´$=\sum_{x\in\X}\,\psi_{\e,x}\,N(x)^{-1}P(x)$\´),\quad and
\item \framebox[10em][l]{\strut$P_\e\neq0\ \lTo\ P|_\X\not\equiv0$}\ .
\end{enumerate}
\end{Satz}\vsp\pagebreak[3]

\begin{Beweis}
In our first proof we use the tensor product property\,\eqref{tens}
and the linearity of the map $P\mto(\Psi(N^{-1}P|_\X))_\e$ to reduce
the problem to the one"=dimensional case. The one"=dimensional case
is covered by the inversion formula\,\ref{sz.itg}\,. Another proof,
following Alon and Tarsi's polynomial method, is described in
Section\,\ref{sec.cn}\´.\pagebreak[2]

Since both sides of the Equation\,\((i)\) are linear in the argument
$P$ it suffices to prove $(X^\d)_\e\,=\,(\Psi(N^{-1}X^\d|_\X))_\e$
in the two cases of Definition\,\ref{def.dl}\,. In each case,
\begin{equation}\ps
\begin{split}
  \Bigl(\Psi(N^{-1}X^\d|_\X)\Bigr)_{\!\e}
  \ceq\biggl(\Psi\Bigl(\bigl(\bigotimes\nolimits_{\!j}N_j^{-1}\bigr)
       \bigotimes\nolimits_{\!j}(X_j^{\d_j}|_{\X_j})\Bigr)\biggr)_{\!\e}\\
  \ceq\biggl(\bigl(\bigotimes\nolimits_{\!j}\Psi^j\,\bigr)
       \bigotimes\nolimits_{\!j}(N_j^{-1}X_j^{\d_j}|_{\X_j})\biggr)_{\!\e}\\
  \ceqby{\eqref{tens}}
       \biggl(\bigotimes\nolimits_{\!j}
              \bigl(\Psi^j(N_j^{-1}X_j^{\d_j}|_{\X_j})\bigr)\,\biggr)_{\!\e}\\
  \ceq\prod_{j\in(n]}\Bigl(\Psi^j(N_j^{-1}X_j^{\d_j}|_{\X_j})\Bigr)_{\!\e_j}\,\ .
\end{split}
\end{equation}
Using the one"=dimensional case of the inversion
formula\,\ref{sz.itg} we also derive
\begin{equation}\ps
  \bigl(\Psi^{j}(N_j^{-1}X_j^{\d_j}|_{\X_j})\bigr)_{\e_j}
  \,=\,(X_j^{\d_j})_{\e_j}
  \,=\,\?{\d_j=\e_j}
  \qquad\text{for all $j\in(n]$ with $\d_j\leq d_j$.}
\end{equation}
Thus in case\,1 ($\fa {j\in(n]}\DP\,\ \d_j=\e_j\leq d_j$\´) ,
\begin{equation}\ps
  \bigl(\Psi(N^{-1}X^{\d}|_\X)\bigr)_\e\,=\,1\,=\,(X^\d)_{\e}\ .
\end{equation}
And in case\,2 ($\ex {j\in(n]}\DP\,\ \e_j\neq\d_j\leq d_j$\´) ,
\begin{equation}\ps
  \bigl(\Psi(N^{-1}X^{\d}|_\X)\bigr)_\e\,=\,0\,=\,(X^\d)_{\e}\ .
  \vspace{-.5\baselineskip}
\end{equation}
\end{Beweis}

Note that the one"=dimensional case of
Theorem\,\ref{sz.cn}\,\((ii)\) is nothing more then the well"=known
fact that polynomials $P(X_1)\neq0$ of degree at most $d_1$ have at
most $d_1$ roots.

With the remark after Definition\,\ref{def.dl}\´, and the knowledge
that $\psi_{d,x}\eqby{\eqref{eqP2}}1$ for all $x\in\X$, we get our
main result as an immediate consequence of Theorem\,\ref{sz.cn}\´:

\index{coefficient formula}
\begin{Satz}[Coefficient formula]\label{sz2.cn}
  Let $\X$ be an integral \(d\)"~grid.
  For each polynomial $P=\sum_{\d\in\N^n}P_\d X^\d\,\in\,\R[X]$
  of total degree $\deg(P)\leq\Sigma d$,
\begin{enumerate}[(i)]
\item \framebox[10em][l]{\strut$P_d\,=\,\Sigma(N^{-1}P|_\X)$}\quad
  (\´$=\sum_{x\in\X}\,N(x)^{-1}P(x)$\´),\quad and
\item \framebox[10em][l]{\strut$P_d\neq0\ \lTo\ P|_\X\not\equiv0$}\ .
\end{enumerate}
\end{Satz}\vspp

\noindent This main theorem looks simpler then the more general
Theorem\,\ref{sz.cn}\´, and you do not have to know the concept of
\(d\)"=leading multiindices to understand it. Furthermore, the
applications in this paper do not really make use of the generality
in Theorem\,\ref{sz.cn}\,. However, we tried to provide as much
generality as possible, and it is of course interesting to
understand the role of the degree restriction in
Theorem\,\ref{sz2.cn}\,.

The most important part of this results, the implication in
Theorem\,\ref{sz2.cn}\,\((ii)\), which is known as Combinatorial
Nullstellensatz was already proven in \cite[Theorem\,1.2]{al2}, for
integral domains. Note that $P_d=0$ whenever $\deg(P)<\Sigma d$, so
that the implication seems to become useless in this situation.
However, one may modify $P\!$, or use smaller sets $\X_j$ (\´and
hence smaller $d_j$\´), and apply the implication then. So, if
$P_\d\neq0$ for a $\d\leq d$ with $\Sigma\d=\deg(P)$ then it still
follows that $P|_\X\not\equiv0$. De facto, such $\d$ are
\(d\)"~leading.

If, on the other hand, $\deg(P)=\Sigma d$, then $P_d$ is, in
general, the only coefficient that allows conclusions on $P|_\X$ as
in Theorem\,\ref{sz2.cn}\,\((ii)\). This follows from the
modification methods of Section\,\ref{sec.cn}\,. More precisely, if
we do not have further information about the \(d\)"~grid $\X$, then
the \(d\)"~leading coefficients are the only coefficients that allow
such conclusions. For special grids $\X$, however, there may be some
other coefficients $P_\d$ with this property, e.g.,\ $P_0$ in the
case $0=(0,\dotsc,0)\in\X$.

Note further that for special grids $\X$, the degree restriction
in Theorem\,\ref{sz2.cn} may be weakened slightly. If, for example,
$\X=\F[q]\n$, then the restriction $\deg(P)\leq\Sigma d+q-2$
suffices; see the 
footnote on page\,\pageref{foo.dr} for an explanation.

The following corollary is a consequence of the simple fact that
vanishing sums \dashed{the case $P_d=0$ in
Theorem\,\ref{sz2.cn}\,(\(i\))} do not have exactly one nonvanishing
summand. It is very useful if a problem possesses exactly one
trivial solution: if we are able to describe the problem by a
polynomial of low degree, we just have to check the degree, and
Corollary\,\ref{kor.cn} guarantees a second (\´in this case,
nontrivial\´) solution. There are many elegant applications of this;
for some examples see Section\,\ref{sec.fap}\,. We will work out a
general working frame in Section\,\ref{sec.asp}\,. We have\´:

\index{theorem!about second nonzeros}
\begin{Korollar}\label{kor.cn}
Let $\X$ be an integral \(d\)"~grid. For polynomials $P$ of degree
$\deg(P)\,<\,\Sigma d$ (or, more generally, for polynomials with
vanishing \(d\)"=leading coefficient $P_d=0$),
%
$$\boxed{\Strut\ \Abs{\{\,x\in\X\mit P(x)\neq0\,\}}\,\neq\,1\ }\ .$$
\end{Korollar}\vspp

\noindent If the grid $\X$ has a special structure -- for example,
if $\X\sb\Rl_{>0}\!\!\strut^n\!$ -- this corollary may also hold for
polynomials $P$ with vanishing \(d\)"=leading coefficient $P_\e=0$
for some $\e\neq d$. The simple idea for the proof of this, which
uses Theorem\,\ref{sz.cn} instead of Theorem\,\ref{sz2.cn}\´, leads
to the modified conclusion that
\begin{equation}
\Abs{\{\,x\in\X\mit \psi_{\e,x}P(x)\neq0\,\}}\,\neq\,1\ .
\end{equation}

Note further that the one"=dimensional case of
Corollary\,\ref{kor.cn} is just a reformulation of the well"=known
fact that polynomials $P(X_1)$ of degree less than $d_1$ do not have
$d_1=\abs{\X_1}-1$ roots, except if $P=0$.

The example $P=2X_1+2\in\Z_4[X_1]$, $\X=\{0,1,-1\}$ shows that
Corollary\,\ref{kor.cn} does not hold over arbitrary grids. However,
if $\X=\Z_m\n=:\R^n$ with $m$ not prime, the grid $\X$ is not
integral; yet assertion\,\ref{kor.cn} holds anyway. Astonishingly,
in this case the degree condition can be dropped, too. We will see
this in Corollary\,\ref{kor3.cn}\,.

We also present another proof of Corollary\,\ref{kor.cn} that uses
only the weaker part\,\((ii)\) of Theorem\,\ref{sz.cn}\,, to
demonstrate that the well"=known Combinatorial Nullstellensatz, our
Theorem\,\ref{sz2.cn}\,\((ii)\),
would suffice for the proof of the main part of the corollary%
\´:
\begin{Beweis}
Suppose $P$ has exactly one nonzero $x_0\in\X$. Then
\begin{equation}\ps
Q\ :=\ P-P(x_0)N^{-1}(x_0)L_{\X,x_0}\ \in\ \R[X]
\end{equation}
vanishes on the whole grid $\X$, but possesses the nonvanishing and
\(d\)"=leading coefficient
\begin{equation}\ps
Q_d\ =\ -P(x_0)N^{-1}(x_0)\ \neq\ 0\,\ ,
\end{equation}
in contradiction to Theorem\,\ref{sz.cn}\,\((ii)\).
%
\end{Beweis}

A further useful corollary, and a version of Chevalley and Warning's
classical result \dashed{Theorem\,\ref{sz.war} in this paper} is the
following result (\´see also \cite{sch2} for a sharpening of
Warning's Theorem, and Theorem\,\ref{kor4.cn} for a similar result
over $\Z_{p^k}\!$\´)\´:\enlargethispage{\baselineskip}

\index{theorem!of Chevalley and Warning}
\begin{Korollar}\label{kor2.cn}
Let $\X\sb\F[p^k]\!\!\n$ be a \(d\)"~grid and
$P_1,\dotsc,P_m\,\in\,\F[p^k][X_1,\dotsc,X_n]$.\\
If $(p^k-1)\sum_{i\in(m]}\deg(P_i)\,<\,\Sigma d$, then
$$
 \boxed{\Strut\
 \Abs{\bigl\{\,x\in\X\,\mit P_1(x)=\dotsb=P_m(x)=0\,\,\bigr\}}\ \neq\ 1\ }\ .
$$
\end{Korollar}\vspp

\begin{Beweis}
Define
\begin{equation}\ps\label{as.cn2}
P\,:=\,\prod_{i\in(m]}(1-P_i^{\,p^k\!-1})\,\ ;
\end{equation}
then for points $x=(x_1,\dotsc,x_n)$\´,
\begin{equation}\ps
 P(x)\neq0\ \ \Longleftrightarrow\ \ \fa i\in(m]\DP\,P_i(x)=0\,\ ,
\end{equation}
and hence
\begin{equation}\ps
 \Abs{\bigl\{\,x\in\X\mit P_1(x)=\dotsb=P_m(x)=0\,\bigr\}}
 \ =\ \Abs{\bigl\{\,x\in\X\mit P(x)\neq0\,\bigr\}}
 \ \eqby[\neq]{\ref{kor.cn}}\ 1\,\ ,
\end{equation}
since
\begin{equation}\ps
  \deg(P)\ \leq\ \sum_{i\in(m]}(p^k-1)\deg(P_i)
         \ <\ \Sigma d\,\ .\vspace{-\baselineskip}
\end{equation}
\end{Beweis}\vspp


\section{First applications and the application principles}\label{sec.fap}

In this section we present some short and elegant examples of how
our theorems may be applied. They are all well"=known, but we wanted
to have some examples to demonstrate the flexibility of these
methods. This flexibility will also be emphasized through the
general working frame described in Section\,\ref{sec.asp}, for which
the applications of this section may serve as examples. Alon used
them already in \cite{al2} to demonstrate the usage of
implication\,\ref{sz2.cn}\,\((ii)\)\,; whereas we prove them by
application of Theorem\,\ref{sz2.cn}\,\((i)\), and corollaries
\ref{kor.cn} and \ref{kor2.cn}\´, an approach which is -- in most
cases -- more straightforward and more elegant. The main advantage
of the coefficient formula\,\ref{sz2.cn}\,\((i)\) can be seen in the
proof of Theorem\,\ref{sz.war}\,, where the
implication\,\ref{sz2.cn}\,\((ii)\) does not suffice to give a proof
of the full theorem. Section\,\ref{sec.mp} will contain another
application that puts the new quantitative aspect of coefficient
formula\,\ref{sz2.cn} into the spotlight.

Our first example was originally proven in \cite{afk}\´:
\index{theorem!about subgraphs}%
\begin{Satz}\label{ex.3r}
Every loopless \(4\)"~regular multigraph plus one edge
$G=(\´V,E\uplus\{e_0\}\´)$ contains a nontrivial \(3\)"~regular
subgraph.
\end{Satz}

See \cite{afk2} and \cite{mozi} for further similar results. The
additional edge $e_0$ in our version is necessary as the example of
a triangle with doubled edges shows.

We give a comprehensive proof in order to outline the principles\´:
\begin{Beweis}
Of course, the empty graph $(\´\o,\o\´)$ is a (\´trivial\´)
\(3\)"~regular subgraph. So there is one ``solution,'' and we just
have to show that there is not exactly one ``solution.'' This is
where Corollary\,\ref{kor2.cn} comes in. Systems of polynomials of
low degree do not have exactly one common zero. Thus, if the
\(3\)"~regular subgraphs correspond to the common zeros of such a
system of polynomials we know that there has to be a second
(\´nontrivial\´) ``solution.''

  The subgraphs without isolated vertices can be identified with the
subsets $S$ of the set of all edges $\Bar E:=E\uplus\{e_0\}$. Now,
an edge $e\in\bar E$ may or may not lie in a subgraph $S\sb\bar E$.
We represent these two possibilities by the numbers $1$ and $0$ in
$\X_e:=\{0,1\}$ (\´the first step in the algebraization\´), we
define
\begin{equation}\label{eq.XS}\ps
 \x(S)
 \,:=\,\bigl(\?{e\in S}\bigr)_{e\in\bar E}
 \ \in\,\X
 \,:=\,\{0,1\}^{\bar E}
 \,\sb\,\F[3]^{\bar E}\,\ .
\end{equation}
With this representation, the subgraphs $S$ correspond to the points
$x=(x_e)$ of the Boolean grid
 $\X:=\{0,1\}^{\bar E}\sb\F[3]^{\bar E}$; and it is easy to see that
the polynomials
\begin{equation}\label{as.3r}\ps
 P_v\,:=\,\sum_{e\ni v}X_e
  \,\in\,\F[3][\,X_e\!\mit\!e\!\in\!\bar E\,]
  \quad\text{for all $v\in V$}
\end{equation}
do the job, i.e., they have sufficient low degrees and the common
zeros $x\in\X$ correspond to the \(3\)"~regular subgraphs. To see
this, we have to check for each vertex $v\in V$ the number
$\Abs{\{\,e\ni v\mit x_e=1\,\}}\leq5$ of edges $e$ connected to $v$
that are ``selected'' by a common zero $x\in\X=\{0,1\}^{\bar E}$\´:
\begin{equation}\ps
 P_v(x)=0
 \quad\Longleftrightarrow\quad
 \sum_{e\ni v}x_e=0
 \quad\Longleftrightarrow\quad
 \Abs{\{\,e\ni v\mit x_e=1\,\}}\in\{0,3\}\,\ .\!
\end{equation}

Furthermore, we have to check the degree condition of
Corollary\,\ref{kor2.cn}\´, and that is where we need the additional
edge $e_0$:
\begin{equation}\ps
 (3^1-1)\sum_{v\in V}\deg(P_v)
 \,=\,2\abs{V}
 \,=\,\abs{E}
 \,<\,\abs{\bar E}
 \,=\,\Sigma d(\X)\ .
\end{equation}
By Corollary\,\ref{kor2.cn}\´, the trivial graph $\o\sb\bar E$
($x=0$) cannot be the only \(3\)"~regular subgraph.
\end{Beweis}

\index{theorem!of Alon and Füredi} The following simple, geometric
result was proven by Alon and Füredi in \cite{alfu}, and answers a
question by Komj\'{a}th. Our proof uses Corollary\,\ref{kor.cn}\´:

\index{theorem!about cube coverings}
\begin{Satz}\label{ex.hp}
Let $H_1,H_2,\dotsc,H_m$ be affine hyperplanes in $\F^n$ (\´$\F$ a
field\´) that cover all vertices of the unit cube $\X:=\{0,1\}^n$
except one, then $m\geq n$.
\end{Satz}

\begin{Beweis}
Let $\sum_{j\in(n]}a_{i,j}X_j=b_i$ be an equation defining $H_i$,
and set
\begin{equation}\label{as.hp}\ps
 P\,:=\,\prod_{i\in(m]}\,\sum_{j\in(n]}(a_{i,j}X_j-b_i)
  \ \in\,\F[][X_1,\dotsc,X_n]\,\ ;
\end{equation}
then for points $x=(x_1,\dotsc,x_n)$;
\begin{equation}\ps
 P(x)\neq0\ \ \Longleftrightarrow\ \
   \bigl(\,\fa i\in(m]\DP\,\sum_{j\in(n]}a_{i,j}x_j\neq b_i\ \bigr)
          \ \ \Longleftrightarrow\ \ x\nin\bigcup_{j\in(m]}H_j\ \ .
\end{equation}
If we now suppose $m<n$, then it follows that
\begin{equation}\ps
\deg(P)\,\leq\,m\,<\,n\,=\,\Sigma d(\X)\,\ ,
\end{equation}
and hence,
\begin{equation}\ps
 \Abs{\X\sm{\textstyle\bigcup_{j\in(m]}}H_j}
 \ =\ \Abs{\{\,x\in\X\mit P(x)\neq0\,\}}
 \ \eqby[\neq]{\ref{kor.cn}}\ 1\,\ .
\end{equation}
This means that there is not one unique uncovered point $x$ in
$\X=\{0,1\}^n$ \ \!-- \ $m<n$ hyperplanes are not enough to achieve
that.
\end{Beweis}

Our next example is a classical result of Chevalley and Warning that
goes back to a conjecture of Dickson and Artin. There are a lot of
different sharpenings to it; see \cite{msck}, \cite{sch2},
Corollary\,\ref{kor2.cn} and Theorem\,\ref{kor4.cn}\,. In the proof
of the classical version, presented below, we do not use the Boolean
grid $\{0,1\}^n\!\!$, as in the last two examples. We also have to
use Theorem\,\ref{sz2.cn}\,\((i)\) instead of its corollaries. What
remains the same as in the proof of the closely related
Corollary\,\ref{kor2.cn}
is that we have to translate a system of equations
into a single inequality
\´:

\index{theorem!of Chevalley and Warning}
\begin{Satz}\label{sz.war}
  Let $p$ be a prime and $P_1,P_2,\dotsc,P_m\,\in\,\F[p^k][X_1,\dotsc,X_n]$.\smallskip

  If $\sum_{i\in(m]}\deg(P_i)<n$, then 
  $$p\Div\Abs{\{\,x\in\F[p^k]\!\!\n\mit P_1(x)=\dotsb=P_m(x)=0\,\}}\,\ ,$$
  and hence the $P_i$ do not have one unique common zero $x$.
\end{Satz}

\begin{Beweis}
  Define
\begin{equation}\ps
P\,:=\,\prod_{i\in(m]}(1-P_i^{p^k-1})\,\ ;
\end{equation}
then
\begin{equation}\ps
  P(x)\,=\,\begin{cases} 1 & \text{if \,$P_1(x)=\dotsb=P_m(x)=0$\,,}\\
                     0 & \text{otherwise}\end{cases}
       \qquad\text{for all $x\in F_{\!p^k}\n\!$,}
\end{equation}
thus, with $\X:=\F[p^k]\!\!\n$,
\begin{equation}\ps
  \Abs{\bigl\{x\in\F[p^k]\!\!\n\mit P_1(x)=\dotsb=P_m(x)=0\bigr\}}\cdot1
  \ =\ \sum_{x\in\X}P(x)
  \ \eqby{\ref{sz2.cn}\atop\ref{eqN}}\ (-1)^n\,(P)_{d(\X)}
  \ \eqby{\eqref{CWdeg}}\ 0\,\ ,
\end{equation}
  where the last two equalities hold as
\begin{equation}\ps\label{CWdeg}
  \deg(P)\ \leq\ (p^k-1)\sum_{i\in(m]}\deg(P_i)
         \ <\ (p^k-1)\,n
         \ =\ \Sigma d(\X)\,\ .\vspace{-\baselineskip}
\end{equation}
\end{Beweis}

The Cauchy-Davenport Theorem is another classical result. It was
first proven by Cauchy in 1813, and has many applications in
additive number theory. The proof of this result is as simple as the
last ones, but here we use the coefficient
formula\,\ref{sz2.cn}\,\((i)\) in the other direction -- we know the
polynomial map $P|_\X$, and use it to determine the coefficient
$P_d$\´:

\index{theorem!of Cauchy and Davenport}
\begin{Satz}\label{ex.cd}
  If $p$ is a prime, and $A$ and $B$ are two nonempty subsets of
  $\Z_p:=\Z/p\Z$, then
$$
\abs{A+B}\geq\min\{\,p\,,\,\abs{A}+\abs{B}-1\,\}\ .
$$
\end{Satz}\vsp

\begin{Beweis}
  We assume $\abs{A+B}\leq\abs{A}+\abs{B}-2$, and must prove $\abs{A+B}\geq p$.
\medskip

  Define
\begin{equation}\ps
P\,:=\prod_{c\in A+B}(X_1+X_2-c)\,\in\,\Z_p[X_1,X_2]\,\ ,
\end{equation}
set
\begin{equation}\ps
\X_1:=A\,\ ,
\end{equation}
and choose a subset
\begin{equation}\ps
\o\neq\X_2\sb B
\end{equation}
of size
\begin{equation}\ps
\abs{\X_2}\,=\,\abs{A+B}-\abs{A}+2\quad(\´\ \leq\,\abs{B}\ )\,\ .
\end{equation}
\smallskip

Now
\begin{equation}\ps
  P|_{\X_1\times\X_2}\,\equiv\,0\,\ ,
\end{equation}
and
\begin{equation}\ps
  \deg(P)\ =\ \abs{A+B}\ =\ \abs{\X_1}+\abs{\X_2}-2\ =\ d_1(\X_1)+d_2(\X_2)\,\ ,
\end{equation}
so that
\begin{equation}\ps
  \tbinom{\abs{A+B}}{d_1}\cdot1
  \ =\ P_{(d_1,\abs{A+B}-d_1)}
  \ =\ P_d
  \ \eqby{\ref{sz2.cn}}\sum_{x\in\X_1\times\X_2\!\!}\!\!\!0
  \ =\ 0
  \ \in\ \Z_p\,\ .
\end{equation}
Hence
\begin{equation}\ps
p\div\tbinom{\abs{A+B}}{d_1}\,\ ,
\end{equation}
and it follows that
\begin{equation}\ps
\abs{A+B}\,\geq\,p\,\ .
\end{equation}
\end{Beweis}

There are some further number"=theoretic applications, for example,
Erd\H{o}s, Ginzburg and Ziv's Theorem, which also can be found in
\cite{al2}.\bigskip


\section{The matrix polynomial -- another application}\label{sec.mp}

In this section we apply our results to the matrix polynomial
$\Pi(AX)$, a generalization of the graph polynomial (see also
\cite{alta2} or \cite{ya}).

We always assume \randd$"A"=(a_{i,j})\,\in\,\R^{m\times n}\!\!$, and
the product of this matrix with the tuple
\randd$"X":=(X_1,\dotsc,X_n)\,\in\,\R[X]\´\,\n$ is
\rand$"AX":=(\sum_{j\in(n]}a_{ij}X_j)_{i\in(m]}$. Now, $\Pi(AX)$ is
defined in accordance with the definition of $\Pi$ in
Section\,\ref{sec.nc}\´, as follows\´:

\index{matrix polynomial}\index{$A\sel{"|\d}$, $\Pi(AX)$, $\per_\d$,
$\pi_{\hspace{-3pt}A}$, $\abs{\sigma^{-1}}$}
\begin{Definition}[Matrix polynomial]\label{def.mp}
The \emph{matrix polynomial} of
$A=(a_{i,j})\,\in\,\R^{m\times n}\!$ is given by
 \rand$$"\Pi(AX)"\ :=\ \prod_{i\in(m]}\,\sum_{j\in(n]}a_{ij}X_j\ \in\,\R[X]\,\ .$$
\end{Definition}\vsp

It turns out that the coefficients of the matrix polynomial are some
kind of permanents\´:

\index{\(\d\)-permanent}\index{permanent}
\begin{Definition}[\(\d\)"~permanent]\label{def.per} For $\d\in\N^n$ we define the
\emph{\(\d\)"~permanent} of $A=(a_{i,j})\,\in\,\R^{m\times n}\!$
through
 \rand$$
 "\per_\d(A)"\ :=\sum_{\sigma\DP(m]\to(n]\!\!\!\atop\abs{\sigma^{-1}}=\d}\!\!%
   \pi_{\!\!A}(\sigma)\,\ ,
$$
where
 \rand\rand$$
 "\pi_{\!\!A}(\sigma)"\,:=\,\prod_{i\in(m]}a_{i,\sigma(i)}
 \quad\ \text{and}\ \quad
 "\abs{\sigma^{-1}}"\,:=\,\bigl(\abs{\sigma^{-1}(j)}\bigr)_{j\in(n]}\ .
$$
\end{Definition}\vsp

Obviously, $\per_\d(A)=0$ if $\Sigma\d\neq m$. If $m=n$ then
\rand$"\per":=\per_{(1,1,\dotsc,1)}$ is the usual permanent; and, if
$\Sigma\d=m$, it is easy to see that
\begin{equation}
\textstyle\bigl(\´\prod_{j\in(n]}\d_j!\´\bigr)\,\per_\d(A)
  \,=\,\per(A\sel{|\d}),
\end{equation}
where \rand$"A\sel{|\d}"$ is a matrix that contains the
$j^\text{th}$ column of $A$ exactly $\d_j$ times. But note that
$\per_\d(A)$ is, in general, not determined by $\per(A\sel{|\d})$.
If, for example, $(\prod_{j\in(n]}\d_j!)\,1=0$ in $\R$, the
\(\d\)"~permanent $\per_\d(A)$ may take arbitrary values, while
$\per(A\sel{|\d})=0$.

As an immediate consequence of the definitions, we have

\begin{Lemma}\label{perd}
$$
 \boxed{\strut\,\Pi(AX)\,=\,\sum\nolimits_{\d\in\N^n}\per_\d(A)\,X^\d}\ .
$$
\end{Lemma}\vsp

The next theorem now easily follows from our main result,
Theorem\,\ref{sz2.cn}\,. It is an integrative generalization of
Alon's \Index{Permanent Lemma} \cite[Section\,8]{al2},
and of Ryser's permanent formula \cite[p.200]{brry},
which follow as the special cases\´:\smallskip\\
--\ \ $m=n$, $d=(1,1,\dots,1)$
      \emph{of the following \ref{sz.pf}\,\((ii)\) over fields},\\
--\ \ $m=n$, $d=(1,1,\dots,1)$, $\X=\{0,1\}^n$, $b=(0,0,\dots,0)$
      \emph{of \ref{sz.pf}\,\((i)\) over fields.}\smallskip\\
We already proved a slightly weaker version for $\X\sb\N^n\sb\R^n$
in \cite[1.14 \& 1.15]{sch}. This proof was based on Ryser's
formula, and is a little more technical.
\cite[1.10]{sch} is the special case $\X=[d]\sb\N^n\sb\R^n\!$, but
you will have to use \ref{eqN}\,\((v)\) to see this. For some
additional tricks over fields of characteristic $p>0$, see
\cite{dev}. We have\´:

\index{permanent formula}
\begin{Satz}[Permanent formula]\label{sz.pf}
Suppose \rand$"A"=(a_{ij})\,\in\,\R^{m\times n}$ and
$b=(b_i)\,\in\,\R^m$ are given, and let $\X\sb\R^n$ be an integral
\(d\)"~grid. If $m\,\leq\,\Sigma d$, then
\begin{enumerate}[(i)]
\item \framebox{\strut$\displaystyle\per_d(A)
   \,=\,\sum\nolimits_{x\in\X}\,N(x)^{-1}\,\Pi(Ax-b)$}\,\ ,\quad
   and
\item \framebox{\Strut$\displaystyle\per_d(A)\,\neq\,0\quad \lTo\quad
 \ex x\in\X\DP\ (Ax)_1\neq b_1\,,\,\dotsc,\,(A x)_m\neq b_m$}\,\ .
\end{enumerate}
\end{Satz}\vsp

\begin{Beweis}
Part\,\((i)\) follows from Theorem\,\ref{sz2.cn}, as
$\deg\bigl(\Pi(AX-b)\bigr)\,=\,m\,\leq\,\Sigma d$, and since
$\bigl(\Pi(AX-b)\bigr)_d\,=\,\bigl(\Pi(AX)\bigr)_d
 \,\eqby{\ref{perd}}\,\per_d(A)$.
Part\,\((ii)\) is a simple consequence of part\,\((i)\).
\end{Beweis}

\index{coloring}\index{\(\d\)-orientation}\index{orientation}We call
an element $x\in\R^m$ with
 $(Ax)_1\neq0$\,, \dots,\,$(A x)_m\,\neq\,0$
a (\´correct\´) \emph{coloring} of $A$,
and a map \rand$"\sigma"\DP(m]\lto(n]$ with
$\pi_{\!\!A}(\sigma)\neq0$ and $\abs{\sigma^{-1}}=\d$ is a
\´\emph{\(\d\)"~orientation} of $A$. 
With this terminology, Theorem\,\ref{sz.pf} describes a connection
between the orientations and the colorings of $A$, and it is not too
difficult to see that this is a sharpening and a generalization of
Alon and Tarsi's Theorem about colorings and orientations of graphs
in \cite{alta}. That is because, in virtue of the embedding
 $\G{\sto}\lmto A(\G{\sto})$ described in \eqref{def.izm} below,
oriented graphs form a subset of the set of matrices, if $-1\neq1$
in $\R$. The resulting sharpening\,\ref{kor.alta} of the Alon"=Tarsi
Theorem contains Scheim's formula for the number of edge
\(r\)-colorings of a planar \(r\)-regular graph as a permanent and
Ellingham and Goddyn's partial solution of the list coloring
conjecture. We briefly elaborate on this; for even more detail, see
\cite{sch}, where we described this for grids
$\X\sb\N^n\sb\R^n\!\!$, and where we pointed out that many other
graph"=theoretic theorems may be formulated for matrices, too.

\index{$A(\G{\sto})$, $\G{\sto}$, $\sTo$, $\sFrom$, $DE_\d$,
$DO_\d$, $EE$, $EO$}Let \randd$"\G{\sto}"=(V,E,\sTo,\sFrom)$ be a
\emph{oriented multigraph} with \emph{vertex set} $V\!$, \emph{edge
set} $E$ and \emph{defining orientations}
\Rand{\({\shortrightarrow}\!\!\!\!\!\hspace{.75pt}{\shortrightarrow},\,
       {\shortleftarrow}\!\!\!\!\!\hspace{.75pt}{\shortleftarrow}\)}%
$\sTo\DP E\lto V$, \randd$e\lmto"\ed{\sto}"$ and $\sFrom\DP
e\lmto\ed{\sfrom}$\!. $\G{\sto}$ shall be \emph{loopless}, so that
$\ed{\sto}\neq\ed{\sfrom}$ for all $e\in E$. We write
 \rand$"v\in e"$ instead of $v\in\{\ed{\sto} ,\ed{\sfrom}\}$
and define the \emph{incidence matrix} \rand$"A(\G{\sto})"$ of
$\G{\sto}$ by
\begin{equation}\label{def.izm}
 A(\G{\sto}):=(a_{e,v})\,\in\,\R^{E\times V}\!,\,\ \text{where}\ \
 a_{e,v}:=\,\?{\ed{\sto}\´\!=\´v}-\,\?{\ed{\sfrom}\´\!=\´v}
 \,\,\in\,\{-1,0,1\}\,\ .
\end{equation}
With this definition, the \emph{orientations} $\sigma\DP E\ni e\lmto
e^\sigma\in e$ and the \emph{colorings} $x\DP V\lto\R$ of $\G{\sto}$
are exactly the orientations and the colorings of $A(\G{\sto})$ as
defined above. The orientations $\sigma$ of $A(\G{\sto})$ have the
special property $\pi_{\!\!A(\G{\sto})}(\sigma)=\pm1$\vspace{1pt}.
According to this, we say that an orientation $\sigma$ of $\G{\sto}$
is \emph{even}/\emph{odd} if $e^\sigma\neq\ed{\sto}$ (\,i.e.,
$e^\sigma=\ed{\sfrom}$) holds for even/odd many edges $e\in E$. We
write \rand$"DE_\d"$/\rand$"DO_\d"$ for the set of even/odd
orientations $\sigma$ of $\G{\sto}$ with
 $\abs{\sigma^{-1}}=\d\in \N^V\!$. With this notation we have\´:

\index{theorem!of Alon and Tarsi}
\begin{Korollar}\label{kor.alta}
Let $\G{\sto}=(V,E,\sTo,\sFrom)$ be a loopless, directed multigraph
and $\X\sb\R^V$ be an integral \(d\)"~grid; where
$d=(d_v)\,\in\,\N^V\!\!$, and $d_v=\abs{\X_v}-1$ for all $v\in
V\!$.\smallskip\\ If $\abs{E}\,\leq\,\Sigma d$, then
\begin{enumerate}[(i)]
\item $\displaystyle\abs{DE_d}-\abs{DO_d}\ =\ \per_d(A(\G{\sto}))
  \ =\ \sum\nolimits_{x\in\X}\,N(x)^{-1}\´
  {\textstyle\prod\nolimits_{e\in
  E}}\bigl(\´x_{\ed{\sto}}-\,x_{\ed{\sfrom}}\bigr)$\,\ ,
\item $\abs{DE_d}\neq\abs{DO_d}\quad\!\lTo\quad\!
 \ex x\in\X\DP\ \!\fa e\in E\DP\
 x_{\ed{\sto}}\neq\,x_{\ed{\sfrom}}$ \ (\´$x$ is a coloring\´).
\end{enumerate}
\end{Korollar}\vsp

Furthermore, it is not so hard to see that, if
\randd$"EE"$/\randd$"EO"$ is the set of even/odd Eulerian subgraphs
of $\G{\sto}$, and $\d:=\abs{\sTo^{-1}}$, we have
$\abs{DE_\d}=\abs{EE}$ and $\abs{DO_\d}=\abs{EO}$; see also
\cite[2.6]{sch}.

Note that even though Corollary\,\ref{kor.alta} looks a little
simpler than \cite[1.14 \& 2.4]{sch}, there is some complexity
hidden in the symbol $N(x)$. If the ``lists'' $\X_v$
(\´$\X=\prod_{v\in V}\X_v$\´) are all equal, this becomes less
complex. Further, if the graph $\G{\sto}$ is the line graph of a
\(r\)"~regular graph, so that its vertex colorings are the edge
colorings of the \(r\)"~regular graph, then the whole right side
becomes very simple. The summands are then \dashed{up to a constant
factor} equal to $\pm1$; or to $0$, if $x=(x_v)_{v\in V}$ is not a
correct coloring. The corres"-pon"-ding specialization of
equation\,\ref{kor.alta}\,\((i)\) was already obtained in
\cite{elgo} and \cite{schm}.

If in addition $\G{\sto}$ is planar this formula becomes even
simpler, so that the whole right side is \dashed{up to a constant
factor} the number of edge \(r\)"~colorings of the \(r\)"~regular
graph. \index{theorem!of Scheim}Scheim\,\cite{schm} proved this
specialization in his approach to the four color problem for
\(3\)"~regular graphs using a result of Vigneron \cite{vig}.
However, with Ellingham and Goddyn's
generalization\,\cite[Theorem\,3.1]{elgo} of Vigneron's result, this
specialization also follows in the \(r\)"~regular
case.\enlargethispage{.4\baselineskip}

As the left side of our equation does not depend on the choice of
the \(d\)"~grid $\X$, the right side does not depend on it, either.
In our special case, where the right side is the number of
\(r\)"~colorings of the line graph of a planar \(r\)"~regular graph,
this means that if there are colorings to equal lists $\X_v$ of size
$r$ (\´e.g., $\X=[r)^V$\´), then there are also colorings to
arbitrary lists $\X_v$ of size $\abs{\X_v}=r$ -- which is just
Ellingham and Goddyn's confirmation of the list coloring conjecture
for planar \(r\)"~regular edge \(r\)"~colorable multigraphs
\cite{elgo}.\bigskip


\section{Algebraically solvable existence problems:\\
         Describing polynomials as equivalent to explicit\\ solutions}\label{sec.asp}
In this section we describe a general working frame to
Theorem\,\ref{sz2.cn}\,\((ii)\) and Corollary\,\ref{kor.cn}, as it
may be used in existence proofs, such as those of 
\ref{kor2.cn}, \ref{ex.hp} or \ref{sz.pf}\,\((ii)\)\,.
We call the polynomials defined in the equations 
\eqref{as.cn2} and \eqref{as.hp} or the matrix polynomial $\Pi(AX)$
in our last example, algebraic solutions, and show that such
algebraic solutions may be seen as equivalent to explicit solutions.
We show that the existence of algebraic solutions, and of nontrivial
explicit solutions are equivalent. To make this more exact, we have
to introduce some definitions. Our definition of problems should not
merely reflect common usage. In fact, the generality gained through
an exaggerated extension of the term ``problem'' through abstraction
is desirable.

\index{problem}\index{solution}\index{$\Pot$, $\S$, $\St$, $\x$}
\begin{Definition}[Problem]\label{def.prob}
A \emph{problem} \rand$"\Pot"$ is a pair \rand$"(\S,\St)"$
consisting of a set $\S$, which we call its set of \emph{solutions};
and a subset $\St\sb\S$, which we call its set of \emph{trivial}
\emph{solutions}.
\end{Definition}

In example\,\ref{ex.3r}\´, the set of solutions $\S$ consists of the
\(3\)"~regular subgraphs and $\St=\{(\o,\o)\}$. These are exact
definitions, but it does not mean that we know if there are
nontrivial solutions, i.e.,  if $\S\neq\St$. The set $\S$ is well
defined, but we do not know what it looks like; indeed, that is the
actual problem.

To apply our theory about polynomials in such general situations, we
have to bring in grids $\X$ in some way. For that, we define
impressions\´:

\index{impression}
\begin{Definition}[Impression]
A triple \rand$"(\R,\X,\x\´)"$ is an \emph{impression} of $\Pot$ if
$\R$ is a commutative ring with $1\neq0$, if
$\X=\X_1\times\dotsb\times\X_n\,\sb\,\R^n$ is a finite integral grid
(\´for some $n\in\N$) and if $\x\DP\S\lto\X$ is a map.
\end{Definition}

As the set $\S$ of solutions is usually unknown, one may ask how the
map $\x\DP\S\lto\X$ can be defined. The answer is that we usually,
define $\x$ on a bigger domain at first, as in
Equation\,\eqref{eq.XS} in example\,\ref{ex.3r}\,.
Then the unknown set of solutions $\S$ (\´more precisely, its image
$\x(\S)$\´) is indirectly described\´:

\index{describing polynomial}
\begin{Definition}[Describing polynomial]
A polynomial $P\in\R[X_1,\dotsc,X_n]$ is a \emph{describing
polynomial} of $\Pot$ over $(\R,\X,\x\´)$ if
$$
\x(\S)\,=\,\supp(P|_\X)\,\ .
$$
\end{Definition}\vsp

\noindent The diagram\,\eqref{sch.sol} in the introduction shows a
schematic illustration of our concept in the case $\St=\o$. The next
question is how it might be possible to reveal the existence of
nontrivial solutions using some knowledge about a describing
polynomial $P\!$, and how to find such an appropriate $P$\!. In view
of our results from Section\,\ref{sec.cf}, we give the following
definition\´:

\index{algebraic solution}
\begin{Definition}[Algebraic solutions]
A describing polynomial $P$ is an \emph{algebraic solution} (\´over
$(\R,\X,\x\´)$\´) of a problem of the form $\Pot=(\S,\o)$ if it
fulfills
$$
deg(P)\,\leq\,\Sigma d(\X)
 \qquad\text{and}\qquad
 P_{d(\X)}\,\neq\,0\,\ .
$$
It is an \emph{algebraic solution} of a problem $\Pot=(\S,\St)$ with
$\St\neq\o$ if it fulfills
$$
\deg(P)\,<\,\Sigma d(\X)
 \ \ \quad\text{and}\!\quad
 \sum_{x\in\x(\St)\!\!}\!\!N(x)^{-1}P(x)\ \neq\ 0\quad(\´\text{e.g., if}\,\
 \abs{\x(\St)}=1\,).
$$
\end{Definition}\vsp

The bad news is that now, we do not have a general recipe for
finding algebraic solutions that indicate the solvability of
problems. However, we have seen that there are several combinatorial
problems that are algebraically solvable in an obvious way. The
construction of algebraic solutions in these examples follows more
or less the same simple pattern, and that constructive approach is
the big advantage. Algebraic solutions are easy to construct if the
problem is not too complex in the sense that the construction does
not require too many multiplications. In many cases algebraic
solutions can be formulated for whole classes of problems, e.g., for
all extended \(4\)"~regular graphs in example\,\ref{ex.3r}, where
the final algebraic solution was hidden in
Corollary\,\ref{kor2.cn}\´.
In these cases a maybe infinite number of algebraic solutions
fit into one single general form, which can be presented on a finite
blackboard or sheet of paper. The concept of algebraic solutions
provides a method of resolution to what we call the ``finite
blackboard problem'', a fundamental problem whenever we view general
situations with infinitely many concrete instances.\pagebreak[3]

If an algebraic solution is found, we can apply
Theorem\,\ref{sz2.cn}\,, Corollary\,\ref{kor.cn} or the following
theorem, which also shows that algebraic solutions always exist,
provided there are nontrivial solutions in the first place.

\begin{Satz}\label{sz.as}
Let $\Pot=(\S,\St)$ be a problem. The following properties are
equivalent\´:
\begin{enumerate}[(i)]
\item There exists a nontrivial solution of $\Pot$; i.e., $\S\neq\St$.
\item 
 There exists an algebraic solution of $\Pot$ over an
 impression\´ \!$(\R,\X,\x\´)$.
\item 
 There exist algebraic solutions of $\Pot$ over each
 impression\´ \!$(\R,\X,\x\´)$ that fulfills either\pagebreak[1]\smallskip\\
   --\ \ $\abs{\R}>2$ \ and \ $\S\neq\St\To\x(\S)\neq\x(\St)$\\
   \phantom{--\ \ $\abs{\R}>2$ \,and \,}(\´e.g., if $\x\!$ is
   injective or if $\St=\o$);
     \quad\ \ \emph{or}\medskip\\
   --\ \ $\abs{\R}=2$ \ and \
   $\abs{\x(\S)}+1\,\equiv\,\abs{\x(\St)}\,\equiv\,\?{\St\neq\o}\pmod2$.
\end{enumerate}
\end{Satz}\vsp

\begin{Beweis}
First, assume \((ii)\), and let $P$ be an algebraic solution. We
want to show that \((i)\) holds. For $\St=\o$, this follows from
Theorem\,\ref{sz2.cn}\,\((ii)\). For $\St\neq\o$, we have
\begin{equation}\ps
  0\ =\ P_d
   \ \eqby{\ref{sz2.cn}}\ \Sigma(N^{-1}P|_\X)
   \ =\sum_{x\in\x(\St)\!\!\!\!}N(x)^{-1}P(x)
      \,\ +\!\!\!\!\sum_{\ \ x\in\x(\S)\sm\x(\St)\!\!\!\!\!\!\!\!\!\!\!\!\!\!\!}
          N(x)^{-1}P(x)\,\ ,
\end{equation}
where the first sum over $\x(\St)$ does not vanish. Hence, the
second sum over the set $\x(\S)\sm\x(\St)$ does not vanish, either.
 Thus $\x(\S)\sm\x(\St)\neq\o$, and $\S\neq\St$ follows.

To prove $(i)\lTo(iii)$ for $\St=\o$, assume \((i)\), and define a
map $y\DP\X\lto\R$ such that $\supp(y)=\x(\S)$ and $\Sigma y\neq0$.
(\´In the case $\abs{\R}=2$, we need $\abs{\x(\S)}\equiv1\pmod2$ to
make this possible.) The interpolation polynomial $P:=(\Psi y)(X)$
to the map $Ny$ described in Theorem\,\ref{sz.ip} now has degree
$deg(P)\leq\Sigma d$, and fulfills
\begin{equation}\ps
\supp(P|_\X)\ \eqby{\ref{sz.ip}}\ \supp(y)
  \ =\ \x(\S)
\end{equation}
and
\begin{equation}\ps
 P_d\ \eqby{\ref{sz2.cn}}\ \Sigma(N^{-1}P|_\X)
    \ \eqby{\ref{sz.ip}}\ \Sigma y
    \ \neq\ 0\ \,.
\end{equation}

To prove $(i)\lTo(iii)$ for $\St\neq\o$, assume \((i)\), and define
a map
 $y\DP\X\lto\R$ such that
 $\supp(y)=\x(\S)$, $\sum_{x\in\x(\St)}y(x)\neq0$ and $\Sigma y=0$.
(\´In the case $\abs{\R}=2$, we need
$\abs{\x(\S)}+1\equiv\abs{\x(\St)}\equiv1\pmod2$ to make this
possible.) Now, the polynomial $P:=(\Psi y)(X)$ has partial degrees
 $\deg_j(P)\leq d_j$, and total degree $deg(P)<\Sigma d$, as
\begin{equation}\ps
 P_d\ \eqby{\ref{sz2.cn}}\ \Sigma(N^{-1}P|_\X)
    \ \eqby{\ref{sz.ip}}\ \Sigma y
    \ =\ 0\,\ .
\end{equation}
It satisfies
\begin{equation}\ps
\supp(P|_\X)\ \eqby{\ref{sz.ip}}\ \supp(y)
  \ =\ \x(\S)
\end{equation}
and
\begin{equation}\ps
 \sum_{x\in\x(\St)\!\!\!\!}\!N(x)^{-1}P(x)
  \ \eqby{\ref{sz.ip}}\sum_{x\in\x(\St)\!\!\!\!}\!y
  \ \neq\ 0\,\ .
\end{equation}\pagebreak[3]

Finally, to show $(iii)\lTo(ii)$, we only have to prove that there
exists an impression $(\R,\X,\x\´)$ as described in \((iii)\). This
is clear, as we may define $\x$ by setting
\begin{equation}\ps
 \x(s):=\begin{cases}
    \,x_{\text{triv}} & \text{for $s\in\St$\,,}\\
    \,x_{\text{good}} & \text{for $s\in\S\sm\St$\,,}
  \end{cases}
\end{equation}
where $x_{\text{triv}}$ and $x_{\text{good}}$ are two distinct,
arbitrary elements in some suitable grid $\X$. 
\end{Beweis}

The arguments in this proof also show that the restrictions to the
impression $(\R,\X,\x\´)$ in part \((iii)\) are really necessary.
If, for example, we had $\abs{\R}=2$, $\St=\o$ and
$\abs{\x(\S)}\,\equiv0\pmod2$, then $P_d=0$, and the problem would
not be algebraically solvable with respect to the impression
$(\R,\X,\x\´)$.\bigskip


\section{The combinatorial nullstellensatz\\
         or how to modify polynomials}\label{sec.cn}

In this section, we describe a sharpening of a specialization of
Hilbert's Nullstellensatz (\´see e.g.\ \cite{dufo}\´), the
so"=called (\´first\´) Combinatorial Nullstellensatz. This theorem,
and the modification methods behind it, can be used in another proof
of the coefficient formulas in Section\,\ref{sec.cf}\,.

We start with an example that illustrates the underlying
\Index{modification method} of this section. It also shows that the
coefficient $P_d$ in Theorem\,\ref{sz2.cn} is, in general, the only
coefficient that is uniquely determined by $P|_\X$\´:

\begin{Beispiel}\label{bsp.cn}
Let $P\in\C[X_1,X_2]$ (\´i.e., $\R:=\C$ and $n:=2$), and define for
$j=1,2$:
\begin{align}
 L_j\!\!\ceq[:=]\!\frac{X_j^5-1}{X_j-1}\,=\,X_j^4+X_j^3+X_j^2+X_j^1+X_j^0
   \qquad\text{\emph{and}}
\\
\X_j\!\!\ceq[:=]\!\{\´x\in\C\mit L_j(x)=0\´\}\,=\,\{x_1,x_2,x_3,x_4\}\ ,
\quad\text{where}\,\ x_k:=e^{\frac{k}{5}2\pi\sqrt{-1}}\ .
\end{align}
Then $d=d(\X)=(3,3)$. Now, for $\e\in\N^2\!$, the polynomial $X^\e
L_1$ (\´and $X^\e L_2$\´) vanishes on $\X$. Therefore, the modified
polynomial
\begin{equation}
P'\,:=\,P+c\,X^\e L_1\,\ ,\ \quad\text{where}\quad c\in\R\ssm0\,\ ,
\end{equation}
fulfills
\begin{equation}
P'|_\X\,=\,P|_\X\,\ ;
\end{equation}
but
the coefficients $P_{\e+(0,0)}$, $P_{\e+(1,0)}$, $P_{\e+(2,0)}$,
$P_{\e+(3,0)}$ and $P_{\e+(4,0)}$ have changed\´:
\begin{equation}
P'_{\e+(i,0)}\,=\,P_{\e+(i,0)}+c\,\neq\,P_{\e+(i,0)}
  \ \quad\text{for}\quad i=0,1,2,3,4\,\ .
\end{equation}
In this way we may modify $P$ without changing the map $P|_\X$.

Now, suppose $\deg(P)\leq\Sigma d=3+3$. Figure\,2 illustrates that
all coefficients $P_\d$ with $\d\leq\Sigma d$ -- except $P_d$ -- can
be modified without losing the condition $\deg P\leq\Sigma d$, so
that they are not uniquely determined by $P|_\X$. If we try to
modify $P_d=P_{(3,3)}$ -- for example, by adding $c\,X^{(0,3)}L_1$
(\´or $c\,X^{(3,0)}L_2$\´) -- we realize that
\begin{equation}
\deg(X^{(0,3)}L_1)\,=\,\deg(X^{(0,7)})\,=\,7\,>\,3+3\,=\,\Sigma d\,\ ,
\end{equation}
and $\deg(P')>\Sigma d$ would follow. The coefficient $P_d$ cannot
be modified in this way.
\begin{figure}
\begin{center}\label{Fig2}
\setlength{\unitlength}{0.00087489in}
\begingroup\makeatletter\ifx\SetFigFont\undefined
\def\x#1#2#3#4#5#6#7\relax{\def\x{#1#2#3#4#5#6}}%
\expandafter\x\fmtname xxxxxx\relax \def\y{splain}%
\ifx\x\y   
\gdef\SetFigFont#1#2#3{%
  \ifnum #1<17\tiny\else \ifnum #1<20\small\else
  \ifnum #1<24\normalsize\else \ifnum #1<29\large\else
  \ifnum #1<34\Large\else \ifnum #1<41\LARGE\else
     \huge\fi\fi\fi\fi\fi\fi
  \csname #3\endcsname}%
\else
\gdef\SetFigFont#1#2#3{\begingroup
  \count@#1\relax \ifnum 25<\count@\count@25\fi
  \def\x{\endgroup\@setsize\SetFigFont{#2pt}}%
  \expandafter\x
    \csname \romannumeral\the\count@ pt\expandafter\endcsname
    \csname @\romannumeral\the\count@ pt\endcsname
  \csname #3\endcsname}%
\fi
\fi\endgroup
{\renewcommand{\dashlinestretch}{30}
\begin{picture}(3000,2800)(0,-150)
\put(2270,390){\shade\ellipse{120}{120}}
\put(1970,390){\shade\ellipse{120}{120}}
\put(1970,690){\shade\ellipse{120}{120}}
\put(1670,390){\shade\ellipse{120}{120}}
\put(1670,690){\shade\ellipse{120}{120}}
\put(1670,990){\shade\ellipse{120}{120}}
\put(1370,390){\blacken\ellipse{120}{120}}
\put(1370,690){\blacken\ellipse{120}{120}}
\put(1370,990){\blacken\ellipse{120}{120}}
\put(1070,390){\blacken\ellipse{120}{120}}
\put(1070,690){\blacken\ellipse{120}{120}}
\put(1070,990){\blacken\ellipse{120}{120}}
\put(1070,1290){\blacken\ellipse{120}{120}}
\put(1070,1590){\shade\ellipse{120}{120}}
\put(770,390){\blacken\ellipse{120}{120}}
\put(770,690){\blacken\ellipse{120}{120}}
\put(770,990){\blacken\ellipse{120}{120}}
\put(770,1290){\blacken\ellipse{120}{120}}
\put(770,1590){\shade\ellipse{120}{120}}
\put(770,1890){\shade\ellipse{120}{120}}
\put(470,390){\blacken\ellipse{120}{120}}
\put(470,690){\blacken\ellipse{120}{120}}
\put(470,990){\blacken\ellipse{120}{120}}
\put(470,1290){\blacken\ellipse{120}{120}}
\put(470,1590){\shade\ellipse{120}{120}}
\put(470,1890){\shade\ellipse{120}{120}}
\put(470,2190){\shade\ellipse{120}{120}}
\path(320,390)(320,2640)
\path(368,2480)(320,2640)(272,2480)
\path(320,2190)(305,2190)
\path(320,1890)(305,1890)
\path(320,1590)(305,1590)
\path(320,1290)(305,1290)
\path(320,990)(305,990)
\path(320,690)(305,690)
\path(320,390)(305,390)
\path(470,240)(2720,240)
\path(2560,192)(2720,240)(2560,288)
\path(2270,240)(2270,225)
\path(1970,240)(1970,225)
\path(1670,240)(1670,225)
\path(1370,240)(1370,225)
\path(1070,240)(1070,225)
\path(770,240)(770,225)
\path(470,240)(470,225)
\thicklines
\put(995,315){\arc{90}{1.5708}{3.1416}}
\put(995,1665){\arc{90}{3.1416}{4.7124}}
\put(1145,1665){\arc{90}{4.7124}{6.2832}}
\put(1145,315){\arc{90}{0}{1.5708}}
\path(950,315)(950,1665)
\path(995,1710)(1145,1710)
\path(1190,1665)(1190,315)
\path(1145,270)(995,270)
\put(395,615){\arc{90}{1.5708}{3.1416}}
\put(395,1965){\arc{90}{3.1416}{4.7124}}
\put(545,1965){\arc{90}{4.7124}{6.2832}}
\put(545,615){\arc{90}{0}{1.5708}}
\path(350,615)(350,1965)
\path(395,2010)(545,2010)
\path(590,1965)(590,615)
\path(545,570)(395,570)
\put(695,615){\arc{90}{1.5708}{3.1416}}
\put(695,765){\arc{90}{3.1416}{4.7124}}
\put(2045,765){\arc{90}{4.7124}{6.2832}}
\put(2045,615){\arc{90}{0}{1.5708}}
\path(650,615)(650,765)
\path(695,810)(2045,810)
\path(2090,765)(2090,615)
\path(2045,570)(695,570)
\put(1370,1290){\ellipse{120}{120}}
\thinlines
\dashline{96.000}(1220,270)(1220,1710)
\dashline{96.000}(350,1140)(1950,1140)
\path(1415,1335)(1527,1447)
\path(2075,795)(2232,952)
\path(1175,1695)(1332,1852)
\path(575,1995)(732,2152)
\put(0,-150){\makebox(0,0)[lb]%
{\smash{{{\SetFigFont{10}{12.0}{sl}Figure\,2\´: Monomials of degree $\leq3+3$.}}}}}
\put(1360,1830){\makebox(0,0)[lb]%
{\smash{{{\SetFigFont{8}{9.6}{rm}$X^{(2,0)}L_2$}}}}}
\put(755,2130){\makebox(0,0)[lb]%
{\smash{{{\SetFigFont{8}{9.6}{rm}$X^{(0,1)}L_2$}}}}}
\put(2255,930){\makebox(0,0)[lb]%
{\smash{{{\SetFigFont{8}{9.6}{rm}$X^{(1,1)}L_1$}}}}}
\put(1550,1425){\makebox(0,0)[lb]%
{\smash{{{\SetFigFont{8}{9.6}{rm}$X^{(3,3)}$}}}}}
\put(-20,2420){\makebox(0,0)[lb]{\smash{{{\SetFigFont{8}{9.6}{rm}$\deg_2$}}}}}
\put(220,2160){\makebox(0,0)[lb]{\smash{{{\SetFigFont{8}{9.6}{rm}6}}}}}
\put(220,1860){\makebox(0,0)[lb]{\smash{{{\SetFigFont{8}{9.6}{rm}5}}}}}
\put(220,1560){\makebox(0,0)[lb]{\smash{{{\SetFigFont{8}{9.6}{rm}4}}}}}
\put(220,1260){\makebox(0,0)[lb]{\smash{{{\SetFigFont{8}{9.6}{rm}3}}}}}
\put(220,960){\makebox(0,0)[lb]{\smash{{{\SetFigFont{8}{9.6}{rm}2}}}}}
\put(220,660){\makebox(0,0)[lb]{\smash{{{\SetFigFont{8}{9.6}{rm}1}}}}}
\put(220,360){\makebox(0,0)[lb]{\smash{{{\SetFigFont{8}{9.6}{rm}0}}}}}
\put(2420,80){\makebox(0,0)[lb]{\smash{{{\SetFigFont{8}{9.6}{rm}$\deg_1$}}}}}
\put(2235,120){\makebox(0,0)[lb]{\smash{{{\SetFigFont{8}{9.6}{rm}6}}}}}
\put(1935,120){\makebox(0,0)[lb]{\smash{{{\SetFigFont{8}{9.6}{rm}5}}}}}
\put(1635,120){\makebox(0,0)[lb]{\smash{{{\SetFigFont{8}{9.6}{rm}4}}}}}
\put(1335,120){\makebox(0,0)[lb]{\smash{{{\SetFigFont{8}{9.6}{rm}3}}}}}
\put(1035,120){\makebox(0,0)[lb]{\smash{{{\SetFigFont{8}{9.6}{rm}2}}}}}
\put(735,120){\makebox(0,0)[lb]{\smash{{{\SetFigFont{8}{9.6}{rm}1}}}}}
\put(435,120){\makebox(0,0)[lb]{\smash{{{\SetFigFont{8}{9.6}{rm}0}}}}}
\end{picture}
}
\end{center}
\end{figure}
\end{Beispiel}

\noindent This example can also be used to illustrate a second proof
of Theorem\,\ref{sz2.cn} (and Theorem\,\ref{sz.cn}\´)\´:

By successive modifications, as above, with
 \rand\begin{equation}\label{eq.Lj}
"L_j"=L_{\X_j}(X_j):=\prod_{\hx\in\X_j}(X_j-\hx)
\end{equation}
\index{trimmed polynomial}\index{$P\!/\´\!\X$, $L_{\X_j}$, $L_j$}in
the general case, it is possible to modify $P$ into a \emph{trimmed}
polynomial \rand$"P\!/\´\!\X"$ with the properties
\begin{equation}\label{eq.trimmed}
P\!/\´\!\X|_\X\,=\,P|_\X\qquad\text{and}\qquad\deg_j(P\!/\´\!\X)\leq
d_j\quad\text{for $j=1,\dotsc,n$.}
\end{equation}
By \ref{eqdef.itg}\((v)\), \(P\!/\´\!\X\ms\) is uniquely determined
if $\X$ is an integral \(d\)"~grid (\´e.g., $\PX[\{x\}]=P(x)$\´). If
$\deg(P)\leq\Sigma d$, then it is obviously possible to leave the
coefficient $P_d$
unchanged during the modification\footnote{\label{foo.dr}%
At this point the degree restriction $\deg(P)\leq\Sigma d$ may be
weakened slightly if the grid $\X$ -- and hence the $L_j$ -- have a
special structure. If, e.g., $L_j=X_j^{k+1}-1$ (\´or
$L_j=X_j^{k+1}-X_j$) for all $j\in(n]$, then $\deg(P)\leq\Sigma
d+k\!\!$ $[\,=(n+1)k\,]$
(\´respectively $\deg(P)\leq\Sigma d+k-1$) suffices.}. %
Therefore we get 
\pagebreak[3]
\begin{equation}\label{eq.Pd}
P_d\,=\,(P\!/\´\!\X)_d
 \,\eqby{\ref{sz.itg}}\,(\Psi(N^{-1}P\!/\´\!\X|_\X))_d
 \,=\,(\Psi(N^{-1}P|_\X))_d
 \,\eqby{\eqref{eqP2}}\,\Sigma(N^{-1}P|_\X)\,\ ;
\end{equation}
and Theorem\,\ref{sz2.cn} follows immediately.

Theorem\,\ref{sz.cn} can also be proven the same way by using the
following obvious generalization (\´Lemma\,\ref{lem.krz}\´) of the
first equation in \eqref{eq.Pd}. Furthermore, we want to mention at
this point that the proof above (\´and the following lemma\´) may
work for some other coefficients $P_\d$ as well if the
$L_j=L_{\X_j}(X_j)$ have a special structure, e.g., $L_j=X^{d_j}-1$.
Of course it works for $P_0$ if $0=(0,\dotsc,0)\in\X$, since all
$L_j$ lack a constant term in this case. Without further information
about the grid $\X$, we ``carry through'' only the \(d\)"~leading
coefficients\´:

\begin{Lemma}\label{lem.krz}
  Let $\X$ 
  be a \(d\)"~grid.
  For each polynomial $P=\sum_{\d\in\N^n}P_\d X^\d\,\in\,\R[X]$ with
  \(d\)"~leading multiindex $\e\leq d\,\in\,\N^n$ (\´e.g., if $\Sigma\e=\deg(P)$\´),
$$
\boxed{\strut\,(P\!/\´\!\X)_\e=P_\e\,}\ .
$$
\end{Lemma}\vspp

If we take a closer look at the modification methods above, we see
that the difference $P-P\!/\´\!\X$ can be written as
\begin{equation}\label{cns}
  P-P\!/\´\!\X\,=\,\sum_{\smash{j\in(n]}}H_jL_j\,\ ,
\end{equation}
with some $H_j\in\R[X]$ of degree $\deg(H_j)\leq\deg(P)-\deg(L_j)$.

If $P|_{\X}\equiv0$, then $P\!/\´\!\X=0$ by the uniqueness of the
trimmed polynomial, and \eqref{cns} yields
$P\,=\,\sum\nolimits_{j\in(n]}H_jL_j$. This was proven for integral
rings in \cite[Theorem\,1.1]{al2}. More formally, we have\´:

\index{Combinatorial Nullstellensatz}
\begin{Satz}[Combinatorial Nullstellensatz]\label{sz3.cn}
  Let $\X=\X_1\times\dotsb\times\X_n\,\sb\,\R^n$ be an integral grid
with associated polynomials
$L_j:=\prod_{\hx\in\X_j}(X_j-\hx)$.\\[2pt]
  For any polynomial $P=\sum_{\d\in\N^n}P_\d X^\d\,\in\,\R[X]$,
  the following are equivalent\´:
\begin{enumerate}[(i)]
\item $\displaystyle P|_{\X}\,\equiv\,0$.
\item $\displaystyle P\!/\´\!\X\,=\,0$.
\item $\displaystyle P\,=\,\sum_{j\in(n]}H_jL_j\,\,\in\,\,\R[X]$ \
  for some polynomials $H_j$ over a ring extension of $\R$.
\item $\displaystyle P\,=\,\sum_{j\in(n]}H_jL_j$ \
  for some $H_j\in\R[X]$ of degree
  $\deg(H_j)\leq\deg(P)-\abs{\X_j}$.
\end{enumerate}
\end{Satz}\vsp

\begin{Beweis}
We already have seen that the implications $(i)\lTo(ii)\lTo(iv)$
hold; and the implications $(iv)\lTo(iii)\lTo(i)$ are trivial.
\end{Beweis}

The implication ``$(i)\lTo(iv)$'' states that polynomials $P$ with
$P|_\X\equiv0$ may be written as $\sum_{j\in(n]}H_jL_j$. In other
words, $P$ lies in the ideal spanned by the polynomials $L_j$. As we
do not know \emph{a priori} that this ideal is a radical ideal,
Hilbert's Nullstellensatz would only provide
$P^k=\sum_{j\in(n]}H_jL_j$ for some $k\geq1$, and without degree
restrictions for the $H_j$ (\´provided $\R$ is an algebraically
closed field\´). Alon suggested calling the stronger (\´with respect
to the special polynomials $L_j$\´) result ``Combinatorial
Nullstellensatz.'' He used it to prove the implication \((ii)\) in
the coefficient formula\,\ref{sz2.cn} \cite[Theorem\,1.2]{al2} and
recycled the phrase ``Combinatorial Nullstellensatz'' for the
implication\,\ref{sz2.cn}\,\((ii)\).\bigskip


\section{Results over $\Z$, $\Z_m$ and other generalizations}\label{sec.ZZm}

There are several ways to generalize the coefficient
formulas\,\ref{sz2.cn} and \ref{sz.cn}\´. This section will address
some of those.

If a grid $\X$ is just affine but we want to use
Theorem\,\ref{sz2.cn}\´, we may apply the homomorphism
\rand$"\pi"\DP r\mto\tfrac{r}{1}$ from $\R$ to the localization
$\R_N$, exactly as in Theorem\,\ref{sz.afg}\´. In particular, this
leads to the implications\´:
\begin{equation}
P_d=1\ \ \lTo\ \ P_d\n[\pi]\neq0\ \ \lTo\ \ P|_\X\not\equiv0\,\ .
\end{equation}

It may also be that there is an integral grid $\hat\X$ over a ring
$\hat\R$, and a homomorphism $\hat\R\lto\R$ that induces a map from
$\hat\X$ into $\X$. Our results may then be applied to a preimage
$\hat P\in\hat\R[X]$ of $P\in\R[X]$. This leads to results about $P$
on not necessarily integral or affine grids $\X$. If, for example,
\rand$\R="\Z_m":=\Z/m\Z$ and $\hat\R=\Z$, we may read the following
formula\,\ref{sz.Zm} modulo $m$ (\´note that it contains only
integer coefficients\´).

\begin{Satz}\label{sz.Zm}
Assume $P\in\Z[X]$ and $\X=[d]:=[d_1]\times\dotsb\times [d_n]$.\\
If $\deg(P)\,\leq\,\Sigma d$, then
$$
 \boxed{\ (-1)^{\Sigma d}\,\bigl[\´{\textstyle\prod_{j\in(n]}}(d_j!)\´\bigr]\,P_d
 \ =\ \sum_{x\in\X}^{\null}\
      \bigl[\´{\textstyle\prod_{j\in(n]}}(-1)^{x_j}\tbinom{d_j}{x_j}\´\bigr]\,P(x)\ }\ .
$$
\end{Satz}\vsp

\begin{Beweis}
This follows from Theorem\,\ref{sz2.cn} and
Lemma\,\ref{eqN}\,\((v)\).
\end{Beweis}

With this theorem we get the following special version of
Corollary\,\ref{kor.cn}\´, which works perfectly well without a
degree condition. (\´See  \cite{must} and \cite{sp} for more
information about polynomial maps $\Z_m\!\n\to\Z_m$.\´)

\pagebreak[4]\index{theorem!about second nonzeros}
\begin{Korollar}\label{kor3.cn}
Let $P\in\Z_m[X]$, and set $\X:=\Z_m\!\n$, 
which we identify with $[m)^n\sb\Z^n\!$.
If $m$ is not prime, and $(m,n)\neq(4,1)$, then\´:
\begin{enumerate}[(i)]
\item \framebox{\Strut$\Abs{\{\,x\in\X\mit P(x)\neq0\,\}}\,\neq\,1$}\ .
\item \framebox{\Strut$P_0\neq0
  \ \,\lTo\ \,\!\ex x\in\X\ssm0\,\DP
          \bigl[\´\prod_{j\in(n]}\binom{m-1}{x_j}\´\bigr]\,P(x)\neq0
  \ \,\lTo\ \,P|_{\X\ssm0}\not\equiv0$}\,.
\item \framebox{\Strut$0
  \,=\,{\displaystyle\sum\nolimits_{x\in\X}}\,
       \bigl[\´\prod\nolimits_{j\in(n]}(-1)^{x_j}\tbinom{m-1}{x_j}\´\bigr]\,P(x)$}\
       .
\end{enumerate}
\end{Korollar}\vspp

\begin{Beweis}
Suppose there is an $\hx\in\X=\Z_m\!\n$ with $P(\hx)\neq0$. By
applying the substitutions $X_j=X_j+\hx_j$, we may assume $0\neq
P(0)=P_0$; and part\,\((i)\) follows from the compounded implication
\((ii)\).\smallskip

Part \((ii)\) follows from part\,\((iii)\), as the summand
$[\prod_{j\in(n]}(-1)^{0}\binom{m-1}{0}]P(0)=P_0$ cannot be the only
nonvanishing summand in the vanishing sum.\smallskip

To prove part\,\((iii)\), we may assume that $P$ has partial degrees
$\deg_j(P)\leq d_j=d_j(\X)$. This is so, as the monic polynomial
$L_j:=\prod_{x\in\X_j}(X_j-x)$ \smallskip maps $\X_j$ to $0$, so
that we may replace $P$ with any polynomial of the form
$P+\sum_{j\in(n]}H_jL_j$ without changing its image $P|_\X$ (\´see
the Example\,\ref{bsp.cn} and \eqref{eq.trimmed} for an illustration
of this method\´). Now let $\hat P\in\Z[X]$ be such that
\begin{equation}\ps
P\,=\,\hat P+m\Z[X]
 \ \,\in\,\,\Z[X]/m\Z[X]
 \,=\,\Z_m[X]
 \qquad\text{and}\qquad
 deg_j(\hat P)\,\leq\,d_j\,\ .
\end{equation}
We only have to show that $m\div(m-1)!^{\,n}\!$, so that the left
side of Equation\,\ref{sz.Zm}\´, applied to $\hat P\!$, vanishes
modulo\! $m$, in the relevant case $d_1,\dotsc,d_n=m-1$:\medskip

If $m\neq4$ and $m=m_1m_2$, with $m_1<m_2<m$, then
$m\div(m-1)!$.\smallskip
%

If $m\neq4$ and $m=p^2$ with $p>2$, then $p<2p<m$ and so
 $m\div p\,(2p)\div(m-1)!$.\smallskip

If $m=4$ and $n\geq2$, then
$m=2^2\div3!^{\,2}\div(m-1)!^{\,n}$.\smallskip

%
\end{Beweis}

The examples $X^3+X+2$ and $X^3-2X^2-X+2\in\Z_4[X]$ show that the
very special case $(m,n)=(4,1)$ in \ref{kor3.cn} is really an
exception. As one can show, these two examples are the only
exceptions to assertion\,(\(i\)) that fulfill the additional
normalization conditions $\deg(P)\leq3$, $P_3\neq-1$ and that the
nonvanishing point is the zero ($P(x)\neq0\,\Eqi\,x=0$\´).\smallskip

We also present another version of Corollary\,\ref{kor2.cn}\,. For
this, we will need the following specialization of
\cite[Lemma\,A.2]{afk2}\´:

\pagebreak[3]\begin{Lemma}\label{lem.p^n}
Let $p\in\N$ be prime, $k>0$ and $c=c(p^k):=\sum_{i\in[k)}(p^i-1)$. 
For $y\in\Z$,
\begin{enumerate}[(i)]
\item $\displaystyle p^c\,\Div\,
       \prod_{0<\hat y<p^k\!\!\!}(y-\hat y)$\ ,\qquad\text{and}
\item $\displaystyle p^{c+1}\,\nDiv\,\prod_{0<\hat y<p^k\!\!\!}(y-\hat y)
       \quad\lEqi\quad p^k\!\div y$\,\ .
\end{enumerate}
\end{Lemma}\vspp

For completeness, we present the relatively short proof\´:
\begin{Beweis}
For each $j\in(k]$ there are exactly $p^{k-j}$ numbers among the
$p^k$ consecutive integers $y$, $y-1$, \dots, $y-(p^k-1)$ that are
dividable by $p^j\!$. Thus\´:\smallskip

If $p^k\!\div y$, then exactly $p^{k-j}-1$ of the factors $y-\hat y$
in the product $\smash{\prod\limits_{0<\hat y<p^k\!\!\!}}(y-\hat y)$
are dividable by $p^j\!$.\smallskip

If $p^k\!\ndiv y$, then at least $p^{k-j}-1$ of these factors are
dividable by $p^j$; and in the case $j=k$, strictly more than
$p^{k-j}-1=0$
are multiples of $p^j=p^k\!$.\\[3pt]
It follows\´:\smallskip

If $p^k\!\div y$, then
 \,$p^c\div\,\prod_{0<\hat y<p^k}(y-\hat y)$,
\,but
 \,$p^{c+1}\ndiv\,\prod_{0<\hat y<p^k}(y-\hat y)$.\smallskip

If $p^k\!\ndiv y$, \,then
 \,$p^{c+1}\div\,\prod_{0<\hat y<p^k}(y-\hat y)$.

\end{Beweis}

The following version of Corollary\,\ref{kor2.cn} (\´see also
Theorem\,\ref{sz.war} and \cite{sch2}\´) reduces to Olson's
Theorem\,\cite[Theorem\,2.1]{afk2}, if
$\deg(P_1)=\dotsb=\deg(P_m)=1$ and if we set $\X:=\{0,1\}^n\!$.
\index{theorem!of Olson}Olson's Theorem can be used, for example, to
prove generalizations of Theorem\,\ref{ex.3r} about regular
subgraphs, such as those in \cite{afk2}. Here we view, more
generally, arbitrary polynomials and arbitrary
\index{grid!\(p\)-integral}\emph{\(p\)"~integral} grids -- i.e.,
grids $\X\sb\Z^n$ with the property\´:
\begin{equation}\label{pr.pint}
\text{For all $j\in(n]$ and all $x,\tx\in\X_j$ with $x\neq\tx\!$,}\quad
  p\ndiv x-\tx\ \,.
\end{equation}
We have\´:

\begin{Satz}\label{kor4.cn}
Let $p\in\N$ be a prime and $\X\sb\Z^n$ a \(p\)"~integral
\(d\)"~grid.\\ For polynomials
$P_1,\dotsc,P_m\,\in\,\Z[X_1,\dotsc,X_n]$, and numbers
$k_1,\dotsc,k_m>0$ small enough so that
$\sum_{i\in(m]}(p^{k_i}-1)\deg(P_i)\,<\,\Sigma d$,
$$
 \boxed{\Strut\
 \Abs{\bigl\{\,x\in\X\,\mit\,\fa i\in(m]\DP\,p^{k_i}\!\div P_i(x)\,\,\bigr\}}\ \neq\ 1\ }\ .
$$
\end{Satz}\vspp

\begin{Beweis}
Set
\begin{equation}\ps
c\,:=\,\sum_{i\in(m]}c_i
 \quad\text{where}\quad
 c_i\,=\,c(p^{k_i})\,:=\,\sum_{j\in[k_i)}(p^j-1)\,\ ,
\end{equation}
define
\begin{equation}\ps
P\,:=\,\prod_{i\in(m]}\,\prod_{0<\hat y<p^{k_i}\!\!\!}(P_i-\hat y)
 \ \in\,\Z[X]
\end{equation}
and let
\begin{equation}\ps
 \bar P\,:=\,P+p^{c+1}\Z[X]\ \in\,\Z[X]/p^{c+1}\Z[X]\,=\,\Z_{p^{c+1}}[X]\,\ .\smallskip
\end{equation}
For points $x=(x_1,\dotsc,x_n)\in\Z^n\!$, set
\begin{equation}\ps
 \bar x\,:=\,(\´x_1+p^{c+1}\Z\,,\dotsc,\,x_n+p^{c+1}\Z\´)\ \in\,(\Z_{p^{c+1}}\!)^n\ ;
\end{equation}
then
\begin{equation}\ps
 \bar\X\,:=\,\{\,\bar x\mit x\in\X\,\}\ \sb\,(\Z_{p^{c+1}}\!)^n
\end{equation}
is an integral \(d\)"~grid, and
 $x\lmto\bar x$ induces a bijection from $\X$ to $\bar\X$.\\[4pt]
Now it follows that
\begin{equation}\ps
\begin{split}
 \bar P(\bar x)\,{\neq}\,0
 \quad\!\ceqby[\Longleftrightarrow]{}\quad\!
 p^{c+1}\ndivps\,P(x)\phantom{\prod_{}}
 \\ \ceqby[\Longleftrightarrow]{\ref{lem.p^n}\´(i)}\quad\!
 \fa i\DP\,p^{c_i+1}\ndivps\,\prod_{0<\hat y<p^{k_i}\!\!}(P_i(x)-\hat y)
 \\ \ceqby[\Longleftrightarrow]{\ref{lem.p^n}\´(ii)}\quad\!
 \fa i\DP\,p^{k_i}\!\div P_i(x)\,\ ,
\end{split}
\end{equation}
and since
\begin{equation}\ps
  \deg(\bar P)\ \leq\ \deg(P)\ \leq\ \sum_{\smash{i\in(m]}}(p^{k_i}-1)\deg(P_i)
         \ <\ \Sigma d\,\ ,
\end{equation}
we obtain
\begin{equation}\ps
 \Abs{\bigl\{\,x\in\X\,\mit\,\fa i\in(m]\DP\,p^{k_i}\!\div P_i(x)\,\,\bigr\}}
 \ =\ \Abs{\bigl\{\,\bar x\in\bar\X\,\mit\,
    \bar P(\bar x)\neq0\,\,\bigr\}}
 \ \eqby[\neq]{\ref{kor.cn}}\ 1\,\ .
\end{equation}

\end{Beweis}

Our result can be generalized further, in the obvious way, by using
\cite[Lemma\,A.2]{afk2}, instead of our Lemma\,\ref{lem.p^n}\´.
However, the result would look a bit more technical. In \cite{afk2}
Alon, Friedland and Kalai also made the following conjecture\´:

\begin{Vermut}
Set $\X:=\{0,1\}^n
$\! and let $P_1,\dotsc,P_m\,\in\,\Z[X_1,\dotsc,X_n]$ be homogenous
polynomials of degree $1$. If $k\in\N$ is small enough so that
$(k-1)\,m\,<\,n$, then
$$
 \boxed{\Strut\
 \Abs{\bigl\{\,x\in\X\,\mit\,\fa i\in(m]\DP\,k\!\div P_i(x)\,\,\bigr\}}\ \neq\ 1\ }\ .
$$
\end{Vermut}\vspp
\pagebreak[3]
\noindent\textbf{Acknowledgement\´:}\\
Many thanks to Alexandra Kallia and Michael Harrison for their help
with the English.


\end{document}